\documentclass{article}
\usepackage{amsfonts,amsmath,amssymb}
\usepackage{graphics, epsfig,enumerate}
\usepackage{psfrag}
\usepackage[page]{appendix}
\usepackage{color}
\allowdisplaybreaks

\newtheorem{theorem}{Theorem}[section]

\newtheorem{lemma}{Lemma}[section]
\newtheorem{corollary}{Corollary}[section]
\newtheorem{remark}{Remark}[section]

\renewcommand{\thesection}{\arabic{section}}
\renewcommand{\theequation}{\thesection.\arabic{equation}}

\numberwithin{equation}{section}
\numberwithin{theorem}{section}
\numberwithin{proposition}{section}
\numberwithin{lemma}{section}
\numberwithin{remark}{section}
\newcommand{\noi}{\noindent}

\newcommand{\dsty}{\displaystyle}

\newcommand{\al}{\alpha}

\newcommand{\be}{\beta}

\newcommand{\gm}{\gamma}
\newcommand{\dl}{\delta}

\newcommand{\lm}{\lambda}

\newcommand{\eps}{\epsilon}
\newcommand{\vp}{\varphi}
\newcommand{\sig}{\sigma}

\newcommand{\om}{\omega}

\newcommand{\nn}{\mathbb{N}}

\newcommand{\rr}{\mathbb{R}}
\newcommand{\rn}{\rr^N}

\newcommand{\bl}[1]{\mathbf{#1}}
\newcommand{\bbox}{\vrule height.6em width.6em 
depth0em} %%%%% Black Box
\newcommand{\dvg}{\operatorname{div}}

\newcommand{\osc}{\operatornamewithlimits{osc}}

\newcommand{\dist}{\operatorname{dist}}

\newcommand{\intl}{\int\limits}

\def\Xint#1{\mathchoice
    {\XXint\displaystyle\textstyle{#1}}%
    {\XXint\textstyle\scriptstyle{#1}}%
    {\XXint\scriptstyle\scriptscriptstyle{#1}}%
    {\XXint\scriptscriptstyle\scriptscriptstyle{#1}}%
    \!\int}
\def\XXint#1#2#3{\setbox0=\hbox{$#1{#2#3}{\int}$}
    \vcenter{\hbox{$#2#3$}}\kern-0.5\wd0}

\def\dashint{\Xint{\raise4pt\hbox to7pt{\hrulefill}}}

\newcommand{\ovl}[3]{\int_{#1}^{#2}\kern-#3pt\raise4pt\hbox to7pt{\hrulefill}\ }

\newcommand{\ovll}[3]{\intl_{#1}^{#2}\kern-#3pt\raise4pt\hbox to7pt{\hrulefill}\ }

\newcommand{\tvl}[2]{\iint_{#1}\kern-#2pt\raise4pt\hbox to7pt{\hrulefill}\ }

\newcommand{\bye}{\end{document}}

\input carleson_25may2012.mac
\input{tcilatex}
\newenvironment{ack}{\medskip{\it Acknowledgement.}}{}
\newcommand*\samethanks[1][\value{footnote}]{\footnotemark[#1]}

\begin{document}

\title{Boundary Estimates for Certain Degenerate\\ and Singular Parabolic
Equations}
\author{Benny Avelin \\
%EndAName
Department of Mathematics, Uppsala Universitet\\
P.O. Box 480, S-75106 Uppsala, Sweden\\
email: \texttt{benny.avelin@math.uu.se} 
\and Ugo Gianazza\thanks{Gianazza and Salsa were partially supported by the 2009 PRIN grant 2009KNZ5FK$\_$002.} \\
%EndAName
Dipartimento di Matematica ``F. Casorati", Universit\`a di Pavia\\
via Ferrata 1, 27100 Pavia, Italy\\
email: \texttt{gianazza@imati.cnr.it} 
\and Sandro Salsa\samethanks  \\
%EndAName
Dipartimento di Matematica ``F. Brioschi"\\
Politecnico di Milano\\
Piazza Leonardo da Vinci 32, 20133 Milano, Italy\\
email: \texttt{sandro.salsa@polimi.it} 
}
\date{\today}

\maketitle

\begin{abstract}
We study the boundary behavior of non-negative solutions 
to a class of degenerate/singular parabolic equations, whose prototype is the parabolic 
$p$--Laplacian. Assuming that such solutions
continuously vanish on some distinguished part of the lateral part $S_T$ of a 
Lipschitz cylinder, we prove Carleson-type estimates, and deduce some 
consequences under additional assumptions on the equation or the domain.
We then prove analogous estimates for non-negative solutions
to a class of degenerate/singular pa\-ra\-bo\-lic equations, of porous medium type.
\vskip.2truecm \noindent\textbf{AMS Subject Classification (2010):} Primary
35K65, 35B65, 35B67; Secondary 35B45 
\vskip.2truecm 
\noindent\textbf{Key
Words:} Degenerate and Singular Parabolic Equations, Harnack estimates,
Boundary Harnack inequality, Carleson estimate.
\end{abstract}

%%%%%%%%%%%%%%%%%%%%%%%%%%%%%%%%%%%%%%%%%%%%%%%%

\vskip.4truecm

%%%%%%%%%%%%%%%%%%%%%%%%%%%%%%%%%%%%%%%%%%%%%%%%%%%

\section{Introduction}\label{S:1}
In this paper we start the study of the boundary behavior of weak solutions to a class of degenerate/singular equations whose prototype is the parabolic $p$--Laplace equation
\begin{equation}\label{Eq:1:1:o}
u_t-\dvg(|D u|^{p-2}D u) =0,  \tag*{(1.1)${}_o$}
\end{equation}
where $Dw$ denotes the gradient of $w$ with respect to the space variables.
Precisely, let $E$ be an open set in $\rn$ and for $T>0$ let $E_T$ denote the
cylindrical domain $E\times(0,T]$. Moreover let
\begin{equation*}
S_T=\partial E\times[0,T],\qquad \partial_P E_T=S_T\cup(E\times\{0\})
\end{equation*}
denote the lateral, and the parabolic boundary respectively.

We shall consider quasi-linear, parabolic partial differential equations of the form
\begin{equation}  \label{Eq:1:1}
u_t-\dvg\bl{A}(x,t,u, Du) = 0\quad \text{ weakly in }\> E_T
\end{equation}
where the function $\bl{A}:E_T\times\rr^{N+1}\to\rn$ is only assumed to be
measurable and subject to the structure conditions
\begin{equation}  \label{Eq:1:2}
\left\{
\begin{array}{l}
\bl{A}(x,t,u,\xi)\cdot \xi\ge C_o|\xi|^p \\
|\bl{A}(x,t,u,\xi)|\le C_1|\xi|^{p-1}%
\end{array}%
\right .\quad \text{ a.e.}\> (x,t)\in E_T,\, \forall\,u\in\rr,\,\forall\xi\in\rn
\end{equation}
where $C_o$ and $C_1$ are given positive constants, and $p>1$.
%\begin{equation}\label{Eq:1:3}
%1<p\le p_*\,\df{=}\,\frac{2N}{N+1}.
%\end{equation}
We refer to the parameters $\datap$ as our structural data, and we write $\gm%
=\gm(N,p,C_o,C_1)$ if $\gm$ can be quantitatively
determined a priori only in terms of the above quantities.
A function
\begin{equation}  \label{Eq:1:4}
u\in C\big([0,T];L^2(E)\big)\cap L^p\big(0,T; W^{1,p}(E)\big)
\end{equation}
is a weak sub(super)-solution to \eqref{Eq:1:1}--\eqref{Eq:1:2} if
for every sub-interval $[t_1,t_2]\subset
(0,T]$
\begin{equation}  \label{Eq:1:5}
\int_E u\vp dx\bigg|_{t_1}^{t_2}+\int_{t_1}^{t_2}\int_E \big[-u\vp_t+\bl{A}%
(x,t,u,Du)\cdot D\vp\big]dxdt\le(\ge)0
\end{equation}
for all non-negative test functions
\begin{equation*}
\vp\in W^{1,2}\big(0,T;L^2(E)\big)\cap L^p\big(0,T;W_o^{1,p}(E)%
\big).
\end{equation*}
This guarantees that all the integrals in \eqref{Eq:1:5} are convergent.

Under the conditions \eqref{Eq:1:2}, equation \eqref{Eq:1:1} is
degenerate when $p>2$ and singular when $1<p<2$, since the modulus of
ellipticity $|Du|^{p-2}$ respectively tends to $0$ or to $+\infty$ as $%
|Du|\to0$. In the latter case, we further distinguish between singular \emph{%
super-critical} range (when $\frac{2N}{N+1}<p<2$), and singular \emph{%
critical and sub-critical} range (when $1<p\le\frac{2N}{N+1}$).
When $p=2$, the equation is uniformly parabolic, and the theory is fairly
complete. %In the
%following particular care will be devoted to the stability
%of the various statements and estimates for $p\approx2$.

For points in $\rn$ we use the notation $\dsty(x^{\prime},x_{\scriptscriptstyle N})$ or
$\dsty(x_1,\dots,x_{\scriptscriptstyle N-1},x_{\scriptscriptstyle N})$; $D'w$
stands for the gradient of $w$ with respect to $x^{\prime}$.

For $y\in \rn$ and $\rho > 0 $, $K_{\rho}(y)$ denotes the cube of edge $2\rho$, centered at $y$
with faces parallel to the coordinate planes. When $y$ is the
origin of $\rn$ we simply write $K_{\rho}$; $%
K^{\prime}_{\rho}(y')$ denotes the $(N-1)$-dimensional cube $\dsty%
\{(x':\,|x_i-y_i|<\rho, i=1,2,...,N-1\}$; we write for short $\dsty%
\{|x_i-y_i|<\rho\}$.

For $\theta > 0$ we also define
\begin{equation*}
Q_{\rho}^{-} (\theta) = K_{\rho} \times (- \theta \rho^p,0], \qquad Q_{\rho}^{+}
(\theta) = K_{\rho} \times (0,\theta \rho^p]
\end{equation*}
and for $(y,s) \in \rn \times \R$
\begin{equation*}
(y,s)+Q_{\rho}^- (\theta) = K_{\rho} (y) \times (s-\theta \rho^p,s], \qquad
(y,s)+Q_{\rho}^+ (\theta) = K_{\rho} (y) \times (s,s+\theta \rho^p].
\end{equation*}

Now fix $(x_o,t_o)\in E_T$ such that $u(x_o,t_o)>0$ and construct the
cylinders
\begin{equation}  \label{Eq:1:6}
(x_o,t_o)+Q_{\rho}^{\pm}(\theta)\quad\text{ where }\quad \theta= \left(\frac{%
c}{u(x_o,t_o)}\right)^{p-2},
\end{equation}
and $c$ is a given positive constant. These cylinders are ``intrinsic'' to
the solution, since their height is determined by the value of $u$ at $%
(x_o,t_o)$. Cylindrical domains of the form $K_{\rho}\times(0,\rho^p]$
reflect the natural, parabolic space-time dilations that leave the
homogeneous, prototype equation \ref{Eq:1:1:o} invariant. The latter however
is not homogeneous with respect to the solution $u$. The time dilation by a
factor $u(x_o,t_o)^{2-p}$ is intended to restore the homogeneity. Most of
the results we describe in this paper hold in such geometry.

Our reference domains are Lipschitz (respectively $C^{1,1}$, $C^2$) domains. We recall that a bounded domain $E \subset \rn$ is said to be a Lipschitz domain, if for each $y\in\partial E$
there exists a radius ${r_o}$, such that in an appropriate coordinate system,
\begin{equation*}%%\label{Eq:1:7}
	\begin{split}
		E\cap K_{8r_o}(y)=\{x=(x',x_{\scriptscriptstyle N})\in\rn : x_{\scriptscriptstyle N} >
	  	\Phi ( x')\}\cap K_{8r_o}(y), \\
	  	\partial E\cap K_{8r_o}(y)=\{x=(x',x_{\scriptscriptstyle N})\in\rn : x_{\scriptscriptstyle N}=
	  	\Phi ( x')\}\cap K_{8r_o}(y),
	\end{split}
\end{equation*}
where $\Phi$ is a Lipschitz function, with $\|D'\Phi\|_{L^\infty}\le L$. 
The quantities $r_o$ and $L$ are independent of $y\in\partial E$. We say
that $L$ is the \emph{Lipschitz constant} of $E$.

\noi We define $C^{1,1}$ and $C^2$ domains analogously, requiring $\Phi$ to be of class
$C^{1,1}$, or $C^2$ respectively. Correspondingly we set $M_{1,1}=\|\Phi\|_{1,1}$,
and $M_2=\|\Phi\|_2$, and we assume they are independent of $y\in \partial E$.
\par
\noi Finally, a Lipschitz (respectively $C^{1,1}$, $C^2$) cylinder 
is a cylindrical domain $E_T$, whose
cross section $E$ is a Lipschitz (respectively $C^{1,1}$, $C^2$) domain.

\smallskip

We are interested in solutions $u$ to \eqref{Eq:1:1}--\eqref{Eq:1:2} continuously \emph{vanishing}
on some distinguished part of the lateral part $S_T$ of a cylinder. Our main goal is to show that near the boundary, $u$ is controlled above in a non-tangential fashion. More precisely,
this means that an inequality of the following type
\begin{equation}\label{Eq:1:8}
u\leq \gamma\, u(P_\rho)
\end{equation}
 holds in a box $\psi_\rho$ of size $\rho$, based on $S_T$, where $P_\rho$ is a non-tangential point and $\gamma$ depends only on the structural data.
The first results of this kind are due to Carleson (\cite{Ca}), for the Laplace
equation in Lipschitz domains and to Kemper (\cite{K}) for the heat equation in domains,
which are locally given by a function satisfying a mixed Lipschitz condition,
with exponent $1$ in the space variables and $\frac12$ in the time variable (also called parabolic 
Lipschitz domains). 
Since then, an inequality like \eqref{Eq:1:8} is known as a Carleson's estimate.

There is another inequality naturally associated to \eqref{Eq:1:8}, namely
\begin{equation} \label{Eq:1:22}
u/v \approx u(P_r)/v(P_r).
\end{equation}
Inequality \eqref{Eq:1:22} is known as the \emph{Boundary Comparison Principle} or the \emph{Boundary Harnack Inequality}. For linear equations, it implies the H\"{o}lder continuity up to the boundary of the quotient $u/v$ and that the vanishing speed of $u$ and $v$ is the same.

Both \eqref{Eq:1:8} and \eqref{Eq:1:22} have been generalized to more general contexts and operators and they have become essential tools in analyzing the boundary behavior of
non-negative solutions.

 In the elliptic context we mention \cite{JK} for the Laplace operator in non-tangentially accessible domains, \cite{CFMS}, and \cite{A},  \cite{B}, \cite{FGMS} for elliptic operators in divergence and non-divergence form, respectively, \cite{LN1}, \cite{LN2}, for the $p$--Laplace operator, \cite{CNP}, \cite{CNP1} for the Kolmogorov operator.

We emphasize that for the Laplacian a Carleson estimate has been proved to be equivalent to the boundary Harnack principle as shown in \cite {Aik}. It would be quite interesting to explore this connection between the two inequalities also in the nonlinear setting.

For parabolic operators, we quote \cite{S1}, \cite{FGS2}, \cite{G}, \cite{FSY} for
cylindrical domains, and \cite{FGS} for parabolic Lipschitz domains.

A classical application of the two inequalities is to Fatou-type theorems, but
even more remarkable, is their role played in the regularity theory of two-phase
free boundary problems, as shown in the two seminal papers \cite{C}, \cite{C1},
where a general strategy to attack the regularity of the free boundary governed by the Laplace operator has been set up.

 This technique has been subsequently extended to stationary problems governed by variable coefficients linear and semilinear operators (\cite{CFS}, \cite{FSa}), to fully nonlinear operators (\cite{F}, \cite{Fe}), and to the $p$--Laplace operator (\cite{LN3}, \cite{LN4}).

The free boundary regularity theory for two-phase parabolic problems is less developed.
For Stefan type problems we mention \cite{CS}, \cite{FSa2}, \cite{DaHa} and the references therein.
In particular \cite{DaHa} deals with a one-phase Stefan problem for the $p$--Laplacian when $p>2$.

Our present paper places itself exactly along this line of research. The Carleson estimate for our singular/degenerate equations is the first piece of information, very useful to analyze the regularity of (e.g.) Lipschitz or flat free boundaries. Thanks to recent development in the field of Harnack inequalities for quasi-linear parabolic equations of $p$--Laplace type (\cite%
{DBGV-acta,DBGV-duke,DBGV-sing07,kuusi}), in Theorems~\ref{Thm:2:1} and \ref%
{Thm:3:1} we extend estimate \eqref{Eq:1:8} to non-negative solutions to
\eqref{Eq:1:1}--\eqref{Eq:1:2} in cylindrical Lipschitz domains. According to the theory developed in the above papers, a Carleson type estimate makes sense only for $p>2N/(N+1)$.

Indeed, in the critical  and  sub-critical range,
explicit counterexamples rule out the possibility of a Harnack inequality. Only so-called
Harnack-type estimates are possible, where, however, the ratio of infimum over supremum 
in proper space-time cylinders depends on the solution itself (for more details, see \cite[Chapter~6, \S~11--15]{DBGV-mono}).  

{The approach developed for linear elliptic equations in \cite{CFMS} and essentially at the same time adapted to the linear parabolic equations in \cite{S1}, to prove the Carleson estimate, is centered around two basic estimates for solutions: the Harnack inequality and the geometric decay of the oscillation of $u$ up to the boundary. To see this, let us consider a non-negative solution in a cylinder, and assume further that the solution vanishes on a part of the lateral boundary, which we assume to be a part of the hyperplane $\{x_{\scriptscriptstyle N}=0\}$. Fix a point $\bar P$ at unit distance from the lateral boundary, by translation we may assume that $\overline P = (e_N,0)$. To simplify even more, 
due to the homogeneity of the equation, without loss of generality we may assume that $u(\overline P) = 1$.
A repeated application of the Harnack inequality in a dyadic fashion gives,
\begin{equation}
	\label{eq:carllinj}
	u(P) \leq H^k u(\overline P) \quad \text{ with } \quad \dist(P,\partial E_T) \geq 2^{-k},
\end{equation}
for $P$ in a boundary space-time box $\Psi^+ = E_T \cap K_1(0) \times (-2,-1)$.
To continue, one defines a sequence of boundary space-time boxes
 $\Psi_k^+$, such that $\Psi_k^+\subset\Psi_{k-1}^+\subset\dots \subset \Psi^+$ for all $k\in\nn$.}
 
{Now, suppose that \eqref{Eq:1:8} does not hold true, i.e. there is no constant $\gm>0$, that depends only on the data,  such that $u(P)\le\gm$ for all $P \in \Psi^+_{1}$. Consequently, there must exist $P_1 \in \Psi^+_1$, such that $u(P_1)>H^{h}$, where $h\in\nn$. 
Inequality \eqref{eq:carllinj} implies that $\dist(P_1,\partial E_T) < 2^{-h}$. If $h$ is chosen large enough, by the geometric decay of the oscillation of $u$ up to the boundary, one deduces the existence of $P_2 \in \Psi^+_{2}$ such that $u(P_2)>H^{h+1}$, and $\dist(P_2,\partial E_T) < 2^{-(h+1)}$. Repeating this yields a sequence of points $\{P_j\}_{j=1}^\infty$ approaching the boundary, whereas the sequence $\{u(P_j)\}_{j=1}^\infty$ blows up: this contradicts the assumption that $u$ vanishes continuously on the boundary, and we conclude
\begin{equation*}
	\sup_{\Psi^+_1} u \leq H^h,
\end{equation*}
which is just \eqref{Eq:1:8} in our setting.}

Although the overall strategy in the nonlinear setting follows the same kind of arguments, its implementation presents a difficulty due to the lack of homogeneity of the equations, and there is also a striking difference between the singular and the degenerate case; this is already reflected in the intrinsic character of the interior Harnack inequality, and it is amplified when approaching the boundary through dyadically shrinking intrinsic cylinders. Concerning the Carleson estimate, its statement in the degenerate case can be considered as the intrinsic version of the analogous statement in the linear uniformly parabolic case. Things are different in the singular super-critical case, where, in general, one can only prove a somewhat weaker estimate (see Theorem~3.1), due to the possibility for a solution to extinguish in finite time. Indeed, the counterexamples in \S~3.2 show that one cannot do any better, unless some control of the interior oscillation of the solution is available (Corollary~3.1).

The difference between the two cases, degenerate and singular super-critical, becomes more evident when one considers the validity of a boundary Harnack principle, even in smooth cylinders. Solutions to the parabolic $p$--Laplace equations can vanish arbitrarily fast in the degenerate case $p>2$ (see \S~2.1), so that no possibility exists to prove a boundary Harnack principle in its generality. On the other hand, in the singular case, the existence of suitable barriers provides a linear behavior. Together with Carleson's estimate, this fact implies almost immediately a Hopf principle and the boundary Harnack inequality.

The last section of this note is devoted to extending all the previous results to non-negative
solutions to a large class of degenerate/singular parabolic equations, whose prototype is the
porous medium equation (see \S~\ref{S:4} for all the details).

In a forthcoming paper we plan to extend the boundary Harnack principle to Lipschitz cylinders.

{\begin{remark}{\normalfont After completing the paper, we learnt that in \cite{KMN}
Kuusi, Mingione and Nystr\"om independently proved a Boundary Harnack inequality, which is similar to the one we give here in Theorem~3.3.}\end{remark}}

\begin{ack} 
{The authors are grateful to the referees for all their comments, 
which helped in improving the paper.}
\end{ack}

\section{The Degenerate Case $p>2$}

\subsection{Main Results}\label{S:2:1}

To describe our main results we need to introduce some further notation. 
Let $E_T$ be a Lipschitz
cylinder and fix $\pto\in S_T$; in a neighbourhood of such a point,
the cross section is represented by the graph 
$\{\dsty(x',x_{\scriptscriptstyle N}):\, x_{\scriptscriptstyle N}=\Phi(x')\}$, 
where $\Phi$ is a Lipschitz function and $\|D'\Phi\|_{\infty}\le L$. 
Without loss of generality, from here on we assume $\Phi(x_o')=0$ and $L\ge1$.

For $\rho\in(0,r_o)$, let $x_\rho=(x_o',2L\rho)$, 
$P_\rho=P_\rho\pto=(x'_{o},2L\rho,t_o)\in E_T$ such that $u\prho>0$. 
Note that $\dist(x_\rho,\partial E)$ is of order $\rho$.
Set
\begin{equation*}
\begin{aligned}
&\Psi^-_{\rho}\pto\\
&=E_T\cap\{(x,t):\,|x_i-x_{o,i}|<{\frac\rho4},\
|x_{\scriptscriptstyle N}|<2L\rho, \ t\in(t_o-\frac{\al+\be}2\theta\rho^p,t_o-\be\theta\rho^p]\}
\end{aligned}
\end{equation*}
where $\dsty\theta=\left[\frac c{u\prho}\right]^{p-2}$, 
with $c$ given in Theorem~\ref{Thm:2:3} below,
and $\al>\be$ are two positive parameters.
We are now ready to state our main result in the degenerate case $p>2$.
\begin{theorem}\label{Thm:2:1}
\emph{(Carleson Estimate, $p>2$)}
Let $u$ be a non-negative,
weak solution to \eqref{Eq:1:1}--\eqref{Eq:1:2} in $E_T$.
Assume that
$$(t_o-\theta(4\rho)^p,t_o+\theta(4\rho)^p]\subset(0,T]$$ and that $u$ vanishes continuously on
$$\partial E\cap \{|x_i-x_{o,i}|<2\rho,\,|x_{\scriptscriptstyle N}|<8L\rho\}\times(t_o-\theta(4\rho)^p,t_o+\theta(4\rho)^p).$$
Then there exist two positive parameters $\al>\be$,
and a constant $\tilde\gm>0$, depending only on
 $p,N,C_o,C_1$ and $L$, such that
\begin{equation}\label{Eq:2:1}
u(x,t)\le\tilde\gm\, u\prho\qquad\text{for every}\  \  (x,t)\in\Psi^-_{\rho}\pto.
\end{equation}
\end{theorem}
{\begin{remark}
{\normalfont Without going too much into details here, let us point out that for the prototype equation \ref{Eq:1:1:o}, estimate\eqref{Eq:2:1}
could be extended from Lipschitz cylinders to a wider class of cylinders $E_T$, whose cross section $E$ is a so-called \emph{N.T.A. domain} (non-tangentially accessible domain). For more particulars, we refer the reader to \cite[\S~12.3]{CS}.}
\end{remark}}
\vskip.2truecm 
\noi Weak solutions to \eqref{Eq:1:1} with zero
Dirichlet boundary conditions on a Lipschitz domain are H\"older continuous
up to the boundary (see, for example, \cite[Chapter~III, Theorem~1.2]{dibe-sv}). 
Combining this result with the previous Carleson estimate,
yields a quantitative estimate on the decay of $u$ at the boundary,
invariant by the intrinsic rescaling
\begin{equation*}
x=x_o+\rho y,\qquad t=t_o+\frac{\rho^p}{u\prho^{p-2}}\tau.
\end{equation*}
\begin{corollary}\label{Cor:2:1}
Under the same assumption of Theorem~\ref{Thm:2:1},
we have
\[
0\le u(x,t)\le\gm\left(\frac{\dist(x,\partial E)}{\rho}\right)^{\mu}\,u\prho,
\]
for every $(x,t)\in\Psi^-_{\frac\rho2}\pto$, where $\mu\in(0,1)$ depends only on 
$p,N,C_o,C_1$ and $L$.
\end{corollary}
\noi If we restrict our attention to solutions to the model equation \ref%
{Eq:1:1:o} and to $C^2$ cylinders, the result of Corollary~\ref{Cor:2:1} can
be strengthened. 
%%%%%%%%%%%%%%%%%%%%%%%%%%%%%%%%%%%%%%%%%%%%%%%%
\begin{theorem}\label{Thm:2:2}
\emph{(Lipschitz Decay)} Let $E_T$ be a $C^2$ cylinder
and $u$ a non-negative,
weak solution to \ref{Eq:1:1:o} in $E_T$.
Let the other assumptions of Theorem~\ref{Thm:2:1} hold.
Then there exist two positive parameters $\al>\be$,
and a constant $\gm>0$, depending only on $p$, $N$, and
the $C^2$--constant $M_{2}$ of $E$,
such that
\begin{equation}\label{Eq:2:2}
0\le\,u(x,t)\le\gm\, \left(\frac{\dist(x,\partial E)}{\rho}\right)u\prho,
\end{equation}
for every 
$$(x,t)\in E_T\cap\left\{|x_i-x_{o,i}|<\frac{\rho}4,\ 
0<x_{\scriptscriptstyle N}< 2M_2\rho\right\}\times\left(t_o-\frac{\al+3\be}4\theta\rho^p,t_o-\be\theta\rho^p\right].$$
\end{theorem}
\vskip.2truecm
\begin{remark}
{\normalfont Estimate \eqref{Eq:2:2} is not  surprising, as it is well-known that
under the assumptions of  Theorem~\ref{Thm:2:2}, solutions are Lipschitz continuous
up to the boundary (see \cite{dibe-sv}). {As a matter of fact, in \cite{lieb1,lieb2} 
Lieberman has proved $C^{1+\al}$ regularity up
to the boundary for solutions of a proper $p$--laplacian type equation, 
with conormal and Dirichlet boundary conditions.}
Relying on the recent papers \cite{Boe,KM1,KM2,KM3},
these results can be extended both to a wider class of degenerate
equations with \emph{differentiable} principal part, which have the same
structure of the $p$--Laplacian and to \emph{less regular}
$C^{1,\al}$ domains. We limited ourselves to (1.1)$_o$ and
$C^2$ domains, mainly to avoid
the introduction of further structural assumptions and technical details. }
\end{remark}
Notice that, in general, the bound below by zero in \eqref{Eq:2:2} cannot be
improved. Indeed, when $p>2$, two explicit solutions to the parabolic $p$%
--Laplacian in the half space $\{x_{\scriptscriptstyle N}\ge0\}$, 
that vanish at $x_{\scriptscriptstyle N}=0$, are
given by
\begin{equation}  \label{Eq:2:3}
\begin{array}{l}
u_1(x,t)=x_{\scriptscriptstyle N}, \\
u_2(x,t)=\left(\frac{p-2}{p^{\frac{p-1}{p-2}}}\right)(T-t)^{-\frac{1}{p-2}%
}x_{\scriptscriptstyle N}^{\frac p{p-2}}.%
\end{array}%
\end{equation}
The power--like behavior, as exhibited
in the second one of \eqref{Eq:2:3}, is not the ``worst'' possible case.
Indeed, let $E=\{-1\le x_i\le1, 0\le x_{\scriptscriptstyle N}\le\frac14\}$, and
consider the following Cauchy-Dirichlet Problem in $E\times[0,T[$:
\begin{equation}  \label{Eq:2:4}
\left\{
\begin{array}{l}
u_t-\dvg(|Du|^{p-2}Du)=0 \\
u(x,0)=C\,T^{-\frac1{p-2}}\exp(-\frac{1}{x_{\scriptscriptstyle N}}) \\
u(x',0,t)=0 \\
u(x',\frac14,t)=C(T-t)^{-\frac1{p-2}}e^{-4} \\
u(x,t)=C(T-t)^{-\frac1{p-2}}\exp(-\frac1{x_{\scriptscriptstyle N}}),\qquad x\in\partial
E\cap\{0<x_{\scriptscriptstyle N}<\frac14\},%
\end{array}
\right.
\end{equation}
where
\begin{equation*}
C=\frac1{2(p-1)(p-2)}\left(\frac{e(p-2)}{2p}\right)^{\frac{2p}{p-2}}.
\end{equation*}
It is easy to check that the function
\begin{equation}  \label{Eq:2:5}
u_3=C(T-t)^{-\frac{1}{p-2}}\exp(-\frac1{x_{\scriptscriptstyle N}}),\qquad x_{\scriptscriptstyle N}>0
\end{equation}
is a super-solution to such a problem. Therefore, the \emph{solution} to the
same problem (which is obviously positive) lies below $u_3$ and approaches
the zero boundary value at $x_{\scriptscriptstyle N}=0$ at least with exponential 
speed. 

Example 
\eqref{Eq:2:4}--\eqref{Eq:2:5} can be further generalized.
Let $\gm\in(0,1)$, $E=\{x_{\scriptscriptstyle N}>0\}$, $T=\frac2\gm-1$:
then
\begin{equation}\label{Eq:2:6}
u(x,t)=\left[\frac{p-2}{p-1}\gm^{\frac1{p-1}}(t+1)
\left(\gm+\frac{x_{\scriptscriptstyle N}-2}{t+1}\right)_+\right]^{\frac{p-1}{p-2}}
\end{equation}
is a solution to \ref{Eq:1:1:o} in $E_T$, and vanishes not only on the boundary
$\{x_{\scriptscriptstyle N}=0\}$, but also in the set $\{0<x_{\scriptscriptstyle N}<2-\gm(t+1),\ \ 0<t<T\}$,
which has positive measure.

\subsection{The Interior Harnack Inequality}\label{SS:2:2}

As we mentioned in the Introduction, our results are strongly based on the
interior Harnack inequalities proved in \cite{DBGV-acta, DBGV-duke,
DBGV-mean, kuusi}, that we recall here. 
\begin{theorem}\label{Thm:2:3}
Let $u$ be a non-negative, weak solution
to \eqref{Eq:1:1}--\eqref{Eq:1:2},
in $E_T$ for $p>2$, $\pto\in E_T$ such that $u\pto>0$.

\noi There exist positive constants
$c$ and $\gm$ depending only on
$p,N,C_o,C_1$, such that for all intrinsic cylinders
$\pto+Q_{2\rho}^{\pm}(\theta)$ as in \eqref{Eq:1:6},
contained in $E_T$,
\begin{equation}\label{Eq:2:7}
\gm^{-1}\sup_{K_{\rho}(x_o)} u(\cdot,t_o-\theta\rho^p) \le
u\pto\le\gm\inf_{K_{\rho}(x_o)} u(\cdot,t_o+\theta\rho^p).
\end{equation}

\noi The constants $\gm$ and $c$ deteriorate as $p\to\infty$ in the sense that
$\gm(p),\ c(p)\to\infty$ as $p\to\infty$; however, they are stable as $p\to2$.
\end{theorem}
\begin{remark}
{\normalfont In all the previously mentioned works, the requirement on the cylinder is that
$\pto+Q_{4\rho}^{\pm}(\theta)\subset E_T$: by a proper adjustment of the parameters
$c$ and $\gm$ we can work under the more restrictive condition we are now assuming.}
\end{remark}

As already pointed out in the Introduction, in \cite{S1}, 
the Carleson estimate is a consequence of the Harnack
inequality of \cite{moser1}, and a
geometric argument, based on the control of the oscillation. 
In particular, a key tool is represented by the so-called
Harnack chain, namely the control on the value of $u(x,t)$
by the value of $u\pto$ with $t<t_o$, thanks to the repeated application of 
the Harnack inequality. 

In \cite{DBGV-duke}, the equivalent statement for solutions
to \eqref{Eq:1:1}--\eqref{Eq:1:2} is given, but a careful examination of the proof shows
that such a result actually holds only for solutions defined in $\rn\times(0,T)$,
and not in a smaller domain $E_T$.
{Although the correct form of 
the Harnack chain for solutions defined in $E_T$, when $E\,{\substack {\not= \\ \subset}}\,\rn$,
can be given, nevertheless, such a result is of no use in the proof of Carleson's estimates,
as there are two different, but equally important obstructions.}

First of all, as \eqref{Eq:2:6} shows too, $u$ can vanish 
and hence prevent any further application 
of the Harnack inequality. 
Indeed, let us consider the following two examples.
\vskip.2truecm
Let $\gm\in(0,1)$; the function
\begin{align*}
u(x,t)=&\left[\frac{p-2}{p-1}\,\gm^{\frac1{p-1}}(t+1)\left(\gm+\frac{x_N-2}{t+1}\right)_+\right]^{\frac{p-1}{p-2}}\\
&+\left[\frac{p-2}{p-1}\,\gm^{\frac1{p-1}}(t+1)\left(\gm-\frac{x_N+2}{t+1}\right)_+\right]^{\frac{p-1}{p-2}}
\end{align*}
is a solution to the parabolic $p$--Laplacian in the set $\rn\times(0,\frac2\gm-1)$ and vanishes in the cone 
\begin{equation*}
\left\{
\begin{aligned}
&0<t<\frac2\gm-1\\
&-(2-\gm(t+1))<x_N<2-\gm(t+1).
\end{aligned}
\right.
\end{equation*}
If we take $(x,t)$ and $\pto$ with $t<t_o$ on opposite sides
of the cone, there is no way to build a Harnack chain that connects the two points.
\vskip.2truecm
Let $\dsty\gm_p=\left(\frac1\lm\right)^{\frac1{p-1}}\frac{p-2}p$, with $\lm=N(p-2)+p$,
consider the cylinder $\{x_{\scriptscriptstyle N}>0\}\times(0,(2\gm_p)^\lm)$ and let $x_1=(0,0,\dots,2)$, $x_2=(0,0,\dots,6)$. The function
\begin{align*}
u(x,t)&=t^{-\frac N\lm}\left[1-\gm_p\left(\frac{|x-x_1|}{t^{\frac1\lm}}\right)^{\frac p{p-1}}\right]_+^{\frac{p-1}{p-2}}\\
&+t^{-\frac N\lm}\left[1-\gm_p\left(\frac{|x-x_2|}{t^{\frac1\lm}}\right)^{\frac p{p-1}}\right]_+^{\frac{p-1}{p-2}}
\end{align*}
is a solution to the parabolic $p$--Laplacian in the indicated cylinder and vanishes 
on its parabolic boundary. Notice that such a solution
is the sum of two Barenblatt functions with poles respectively at $x_1$ and $x_2$ and masses $M_1=M_2=1$: in the interval
$0<t<(2\gm_p)^\lm$ the support of $u$ is given by two disjoint regions $R_1$ and $R_2$, and
only at time $T=(2\gm_p)^\lm$ the support of $u$ finally becomes a simply connected set. Once more,
taking $(x,t)$ and $\pto$ respectively in $R_1$ and $R_2$, there is no way to connect them
with a Harnack chain. As a matter of fact, before the two supports touch, each Barenblatt function  does not \emph{feel} in any way the presence of the other one.
In particular, we can change the mass of the two Barenblatt functions: this will modify the time $T$
the two supports touch, but up to $T$, there is no way one Barenblatt component can detect the change performed on the other one.  
\vskip.2truecm
On the other hand, one could think that if we have a solution vanishing on a flat piece of the boundary and strictly positive everywhere in the interior, then one could build a Harnack chain extending arbitrarily close to the boundary. However, this is not the case, as clearly shown
by the following example.

Let us consider a domain $E\subset \rn$, which has a part of its boundary which coincides with the hyperplane $\{x_{\scriptscriptstyle N} = 0 \}$, and let $\Gamma = \partial E \cap \{x_{\scriptscriptstyle N} = 0 \}$. Let $\bar T > 0$, be given and consider a non-negative solution $u$ to
\begin{equation*}
\left\{
\begin{array}{l}
	u_t - \dvg (|D u|^{p-2} Du) = 0, \quad \text{in $E_{\bar T}$} \\
	u > 0, \quad \text{in $E_{\bar T}$} \\
	u = 0, \quad \text{on $\Gamma\times(0,\bar T]$}.
\end{array}
\right.
\end{equation*}
Let $u$ be such that its value is bounded above by the distance to the flat boundary piece 
raised to some given power $a > 0$, i.e.
\begin{equation}\label{Eq:2:8}
	u(x,t) \leq \gm\, \dist(x,\Gamma)^a, \quad a > 0,\ (x,t) \in E_{\bar T},
\end{equation}
where $\gm>0$ is a proper parameter.
Indeed, this is the case of solutions given in \eqref{Eq:2:3} and \eqref{Eq:2:5} with $\bar T=\frac T2$ for example,
and therefore such a situation does take place.

Let $\pto = (x_o',x_{o,{\scriptscriptstyle N}},t_o) \in E_T$ be such that $\dist(x_o, \Gamma) = 1$. The goal is to form a Harnack chain of dyadic non-tangential cylinders approaching the boundary, while the chain stays inside $E_{\bar T}$: we want to control the size of the time interval, which we need to span in order to complete the chain.
Let 
\begin{eqnarray*}
	u_o &=& u\pto \\
	r_k &=& 2^{-k} \\
	x_k &=& (\hat x_o', 2^{-k}) \\
	t_k &=& t_o - c^{p-2} \sum_{i=0}^{k-1} u_i^{2-p} r_i^p\\
	u_k &=& u(x_k,t_k) \approx (2^{-k})^a
\end{eqnarray*}
for $k = 1,\ldots$
Assuming that at each step one can use Harnack's inequality, we get an estimate on the size of $t_k$ from above
\begin{eqnarray*}
	t_k \leq t_o - c^{p-2} \sum_{i=0}^{k-1} (2^{-ai})^{2-p} 2^{-ip}
	\leq t_o - c^{p-2} \sum_{i=0}^{k-1} 2^{ai(p-2)-ip}
\end{eqnarray*}
which diverges to $- \infty$ as $k\to\infty$ and $x_k\to\Gamma$,
if $a \geq p/(p-2)$.
Considering the solution $u_2$ from \eqref{Eq:2:3},
we see that the above dyadic Harnack chain would diverge for such a  solution 
as $a = \frac{p}{p-2}$.

\noi Notice that this is a counterexample to the use of the Harnack chain in the \emph{proof}
of the Carleson estimate, but not a counterexample to the Carleson estimate itself. 

\begin{remark}
{\normalfont The infinite length of the time interval needed to reach the boundary,
is just one face (i. e. consequence) of the finite speed of propagation when $p>2$.
Points $(x,t)$ that lie \emph{inside} a proper $p$--paraboloid 
centered at $\pto$ can be reached, starting from $\pto$: if $u_o$ is very small, 
and therefore the $p$--paraboloid is very narrow,
with small values of $r$ one ends up with very large values of $t$.
On the other hand, points $(x,t)$ that lie \emph{outside} the same $p$--paraboloid 
centered at $\pto$ cannot be reached%: we will comment on  this in Appendix~C
.}
\end{remark}

\begin{remark}
{\normalfont We conjecture that $\frac{p}{p-2}$ is a sort of threshold exponent: when in \eqref{Eq:2:8}
$a<\frac{p}{p-2}$, the regularizing effects of the diffusion kicks in, and eventually the solution
becomes linear, allowing for more precise bounds from below in \eqref{Eq:2:2}; whenever 
$a$ is larger, the time evolution part wins, and we have much less regularity. 
We probably need further information on the behavior of $u$ on the lower
base of the cylinder, in order to make all the previous heuristics more rigorous.
This will be the object of future investigation.}
\end{remark}

Eventually, as stated in Theorem~\ref{Thm:2:1}, the Carleson estimate does hold true.
However, as we could not make use of any form of Harnack chain in its proof, 
we resorted to a contradiction argument, detailed in \S~\ref{SS:2:4:2}
%%%%%%%%%%%%%%%%%%%%%%%%%%%%%%%%%%%%%%%%%%%%%%%%

\subsection{H\"older Continuity and Oscillation Control}

In this subsection we consider $p>1$, since the statements 
are the same in both cases $p>2$ and $1<p<2$. It is well known 
that locally bounded, weak solutions to \eqref{Eq:1:1}--\eqref{Eq:1:2} 
are locally H\"older continuous.
For the full statement and the proof, see \cite[Chapters~III and IV]{dibe-sv}, 
or \cite{GSV} (when $p>2$), for a simpler approach,
based on the same ideas used in \cite{DBGV-acta}. The local H\"older
continuity is a consequence of the following Lemma. Note the slight difference 
between the case $p>2$ and $1<p<2$.

In the following let $u$ be a weak solution to \eqref{Eq:1:1}-\eqref{Eq:1:2}, in $E_T$ for $1 < p < \infty$.
Fix a point in $E_{T}$, which, up to a translation, we may take to be the
origin of $\rr^{N+1}$. For $\rho >0$ consider the cylinder
\begin{equation*}
Q_{\ast }=\left\{
\begin{array}{l}
K_{\rho }\times (-\rho ^{2},0]\quad\text{ if }\, p>2,\\
K_{\rho }\times (-\rho ^{p},0]\quad\text{ if }\, 1<p<2,
\end{array}%
\right.
\end{equation*}%
with vertex at $(0,0)$, and set
\begin{equation*}
\mu _{o}^{+}=\sup_{Q_{\ast }}u,\qquad \mu _{o}^{-}=\inf_{Q_{\ast }}u,\qquad %
\om_{o}=\osc_{Q_{\ast }}u=\mu _{o}^{+}-\mu _{o}^{-}.
\end{equation*}%
 Relying on $\om_{o}$, construct the cylinder
\begin{equation*}
Q_{o}=K_{\rho }\times (-\theta _{o}\rho ^{p},0],\quad \hbox{\rm where}\
\theta _{o}=\left( \frac{c}{\om_{o}}\right) ^{p-2}.
\end{equation*}%
If $p>2$, $c$ being the constant that appears in \eqref{Eq:1:6}, we assume that
\begin{equation}
\om_{o}>c\rho .  \label{Eq:2:13}
\end{equation}
If $1<p<2$, we assume that $\om_o\le1$.
 The previous two condition ensures that $Q_{o}\subset Q_{\ast }$, and the following lemma holds.
 
\begin{lemma}\label{Lm:2:3}
\emph{(H\"older continuity, \cite{DBGV-acta})} There exist constants $\eps$,
$\delta\in(0,1)$, and $c\ge1$, %if $1<p<2$,
depending only on $p,N,C_o,C_1$ such that, setting
\begin{equation*}
\om_n = \dl\om_{n-1}, \quad
\theta_n = \left ( \frac{c}{\om_n} \right )^{p-2}, \quad
\rho_n = \eps \rho_{n-1} \quad \textrm{and} \quad
Q_n = Q_{\rho_n}^- (\theta_n),
\end{equation*}
for all non-negative integers $n$, there holds $Q_{n+1}\subset Q_n$ and
\begin{equation}\label{Eq:2:14}
\osc\limits_{Q_n} u \le \om_n.
\end{equation}
\end{lemma}

\begin{remark}
{\normalfont Given $u$ and $p > 2$, we can directly assume \eqref{Eq:2:13}, since otherwise there is nothing to prove; {moreover, in this case, we have $c=1$. If $1 < p < 2$, $c$ is a suitable constant greater than $1$.}
%: then, in this case, since $u$ is bounded, we can directly assume that $\omega_0 \leq 1$.
Also, note that \eqref{Eq:2:14} yields
\begin{equation*}
\osc_{Q_n} u\le\dl^n\osc_{Q_*} u.
\end{equation*}
Therefore, the lemma builds a sequence of intrinsic cylinders, where the oscillation is controlled
by a proper power of an absolute constant; the starting cylinder is the only one that
does not have an intrinsic size.}
\end{remark}

In the proof of Theorems~\ref{Thm:2:1}--\ref{Thm:2:2} we need the following
two Lemmas. {The former is the well-known reflection principle, whose proof 
is standard. The interested reader can refer to \cite[Lemma~2.7]{LunNys},
or to \cite[Lemma~2.8]{Boe}}. The latter is an alternative DeGiorgi-type Lemma
with initial data, taken from \cite{DBGV-mono}, 
to which we refer for the proof.
\begin{lemma}\label{Lm:2:4}
Let
\begin{align*}
&Q_p = \left\{ x : |x_i| < 1,\ 0 < x_{\scriptscriptstyle N} < 2,\ t_1 < t < t_2 \right\}\\
&Q_n = \left\{ x : |x_i| < 1, -2 < x_{\scriptscriptstyle N} < 0 ,\ t_1 < t < t_2 \right\},
\end{align*}
and $u$ be a non-negative, weak solution to \eqref{Eq:1:1}--\eqref{Eq:1:2},
in $Q_p$, with 	$u = 0$ on $\partial Q_p \cap \{x:\,x_{\scriptscriptstyle N} = 0\}$.
Let
\[
Q' = \left\{ x : |x_i| < 1,\ |x_{\scriptscriptstyle N}| < 2,\ t_1 < t < t_2 \right\},
\]
and define $\tilde u$, $\tilde{\bl{A}}$ as
\begin{align*}%%\label{Eq:2:15}
&\tilde u(x',x_{\scriptscriptstyle N},t)=
\begin{cases}
u(x',x_{\scriptscriptstyle N},t) & \textrm{ if $x_{\scriptscriptstyle N} \geq 0$ } \\
-u(x',-x_{\scriptscriptstyle N},t) & \textrm{ if $x_{\scriptscriptstyle N} < 0$ },
\end{cases}\nonumber\\
%\end{equation}
%\begin{equation}\label{Eq:2:16}
&\tilde{A}_i(x',x_{\scriptscriptstyle N},t)=
\begin{cases}
A_i(x',x_{\scriptscriptstyle N},t) &  \textrm{ if $x_{\scriptscriptstyle N} \geq 0$ }\\
-A_i(x',-x_{\scriptscriptstyle N},t) &  \textrm{ if $x_{\scriptscriptstyle N}<0$ },
\end{cases}
\qquad i=1,\dots,N-1,\\
%\end{equation}
%
%\begin{equation}\label{Eq:2:17}
&\tilde{A}_{\scriptscriptstyle N}(x',x_{\scriptscriptstyle N},t)=
\begin{cases}
A_{\scriptscriptstyle N}(x',x_{\scriptscriptstyle N},t) &  \textrm{ if $x_{\scriptscriptstyle N} \geq 0$ }\\
A_{\scriptscriptstyle N}(x',-x_{\scriptscriptstyle N},t) &  \textrm{ if $x_{\scriptscriptstyle N}<0$ },
\end{cases}\nonumber
\end{align*}
where $\tilde{\bl{A}}(x',x_{\scriptscriptstyle N},t)=\tilde{\bl{A}}(x',x_{\scriptscriptstyle N},t,
\tilde u(x',x_{\scriptscriptstyle N},t),D\tilde u(x',x_{\scriptscriptstyle N},t))$.
Then $\tilde{u}$ is a weak solution in $Q'$ to \eqref{Eq:1:1}--\eqref{Eq:1:2},
with $\bl{A}$ substituted by $\tilde{\bl{A}}$.
\end{lemma}

\begin{lemma}\label{Lm:2:4bis}
Let $u$ be a {\it non-negative}, weak super-solution to (\ref{Eq:1:1})--(\ref{Eq:1:2}) 
in $E_T$, and let $a \in (0,1)$ be given. Let $Q^+=K_{2\rho}(y)\times(s,s+\theta(2\rho)^p]$, and 
$\xi$ be  a positive number such that 
\begin{equation*}
u(x,s)\ge \xi\quad\text{ for a.e. }\> x\in K_{2\rho}(y)
\end{equation*}
and 
\begin{equation}\label{Eq:2:15}
\frac{|[u<\xi]\cap Q^+|}{|Q^+|}
\le\dl\frac{\xi^{2-p}}{\theta}
\end{equation}
for a constant $\dl\in(0,1)$ depending only on the data 
and $a$, and independent of $\xi$, $\rho$, and 
$\theta$.
Then 
\begin{equation*}%%\label{Eq:1:4:4}
u\ge a\xi\quad\text{ a.e. in }\>  K_{\rho}(y)\times(s,s+\theta(2\rho)^p].
\end{equation*}
\end{lemma}
Notice that \eqref{Eq:2:15} is automatically satisfied, by taking $\theta=\frac\dl{\xi^{p-2}}$.
%%%%%%%%%%%%%%%%%%%%%%%%%%%%%%%%%%%%%%%%%%%%%%%
\subsection{Proofs of the Theorems~\protect\ref{Thm:2:1}--\protect\ref%
{Thm:2:2}}

\subsubsection{\textit{Flattening the Boundary}}\label{S:2:4:1} 
If we introduce the new variables
\begin{equation*}
y_i=x_i,\quad i=1,\dots,N-1,\qquad y_{\scriptscriptstyle N}=x_{\scriptscriptstyle N}-\Phi(x^{\prime}),
\end{equation*}
the portion $\partial E\cap \{|x_i-x_{o,i}|<2\rho,\ |x_{\scriptscriptstyle N}|<8L\rho\}$ 
coincides with a portion of the hyperplane $y_{\scriptscriptstyle N}=0$. 
Let $\tilde{K}_{2\rho}(x_o)=\{|x_i-x_{o,i}|<2\rho,\ |x_{\scriptscriptstyle N}|<4L\rho\}$.
We orient $y_{\scriptscriptstyle N}$ so that 
$E\cap \tilde{K}_{2\rho}(x_o) \subset\{y_{\scriptscriptstyle N}>0\}$. It is easy to see that,
with respect to
the new variables, \eqref{Eq:1:1} becomes
\begin{equation*}
u_t-\dvg_y\tilde{\bl{A}}(y,t,u,D_yu)=0,
\end{equation*}
and $\dsty\tilde{\bl{A}}(y,t,u,D_yu)$ satisfies the same kind of structural conditions as
given in \eqref{Eq:1:2}. {We refer for more details to \cite[Chapter~X, \S~2]{dibe-sv}}.
%All the details are in \S~B of the Appendix.

Denoting again by $x$ the transformed variables $y$, under the previously
described change of variables, $u$ is still a solution to an equation of
type \eqref{Eq:1:1}--\eqref{Eq:1:2}, and satisfies a homogeneous Dirichlet condition
on a \emph{flat} boundary. This is a notion, which might not 
be standard in the literature,
and we introduce it, mainly to simplify the notation in the following:
consider the set $\tilde{K}_{2\rho}(x_o)$, and the set
\begin{equation}\label{Eq:2:16}
K^*_{2\rho}(x_o)=\{|x_i-x_{o,i}|<2\rho, 0<x_{\scriptscriptstyle N}<4L\rho\}.
\end{equation}

\begin{definition}%%\label{Def:2:1}
We say that the boundary of $E$ is {\bf flat with respect to $x_{\scriptscriptstyle N}$} if
the portion $\partial E\cap {\tilde{K}_{2\rho}(x_o)}$ coincides with the portion
of the hyperplane $\{x_{\scriptscriptstyle N}=0\}\cap{\tilde{K}_{2\rho}(x_o)}$, and $K^*_{2\rho}(x_o)\subset E$.
We orient $x_{\scriptscriptstyle N}$ so that $E\cap \tilde{K}_{2\rho}(x_o)\equiv K^*_{2\rho}(x_o)\subset\{x_{\scriptscriptstyle N}>0\}$.
\end{definition}

\noi Therefore, proving Theorem~\ref{Thm:2:1} reduces 
to the proof of the following lemma.
\begin{lemma}
Let $u$ be a non-negative, weak solution to \eqref{Eq:1:1}--\eqref{Eq:1:2} in $E_T$
for $p>2$. Take $\pto\in S_T$, $\rho\in(0,r_o)$, let $P_\rho=(x_o',2L\rho,t_o)$,
and assume that $u\prho>0$, $\partial E$ is flat
with respect to $x_{\scriptscriptstyle N}$, and $(t_o-\theta(4\rho)^p,t_o+\theta(4\rho)^p]\subset(0,T]$,
where $\theta=\left[\frac c{u\prho}\right]^{p-2}$, with $c$ given by Theorem~\ref{Thm:2:3}.
Suppose that $u$ vanishes continuously on $(\partial E\cap{K_{2\rho}(x_o)})\times(t_o-\theta(4\rho)^p,t_o+\theta(4\rho)^p]$.

Then there exist two positive parameters $\al>\be$, and a constant $\tilde\gm>0$ 
depending only on $\datap$,
such that
\begin{equation}\label{Eq:2:18}
u(x,t)\le\tilde\gm\, u\prho
\end{equation}
$$\forall(x,t)\in\left\{|x_i-x_{o,i}|<\frac\rho4, \
0<x_{\scriptscriptstyle N}<2L\rho\right\}\times\left(t_o-\frac{\al+\be}2\theta\rho^p,t_o-\be\theta\rho^p\right].$$
\end{lemma}

\subsubsection{\textit{Intrinsic Rescaling and Harnack-based Upper Bounds}}\label{SS:2:4:2}

The change of variable
\begin{equation*}
x\,\to\,\frac{x-x_o}{2L\rho},
\qquad t\,\to\,u\prho^{p-2}\frac{t-t_o}{\rho^p}
\end{equation*}
maps $Q^*_{\rho}(4^p\theta)=K^*_{2\rho}(x_o)\times(t_o-\theta(4\rho)^p,t_o+%
\theta(4\rho)^p]$ into
\begin{equation*}
\hat Q=\{|y_i|<\frac1L
,\ 0<y_{\scriptscriptstyle N}<2\}\times(-4^pc^{p-2},4^pc^{p-2}],
\end{equation*}
$x_\rho$ into $y_o=(0,\dots,0,1)$, $\tilde{K}_{\rho}(x_o)$ into $\tilde K_1=\{|y_i|<\frac1{2L}, |y_{\scriptscriptstyle N}|<1\}$, 
$K^*_{2\rho}(x_o)$ into $K^*_2(y_o)=\{|y_i|<\frac1L,\ 0<y_{\scriptscriptstyle N}<2\}$ and the
portion of the lateral boundary $S_T\cap Q^*_{\rho}(4^p\theta)$ into
\begin{equation*}
 \Xi=\left\{(y^{\prime},0):\ |y_i|<\frac1L\right\}\times\left(-4^pc^{p-2},4^pc^{p-2}\right].
\end{equation*}
Denoting again by $(x,t)$ the transformed variables and letting $y_o=(0,\dots,0,1)$, 
the rescaled function
\begin{equation*}
v_\rho(x,t)=\frac{1}{u\prho}u(2L\rho\, x+x_o,t_o+\frac{t\rho^p}{u\prho^{p-2}})
\end{equation*}
is a non-negative, weak solution to
\begin{equation*}
\partial_t v_\rho-\dvg{\bl{A}}_\rho(x,t,v,Dv)=0
\end{equation*}
in $\hat Q$, where $v_\rho(y_o,0)=1$, and it is easy to see that ${\bl{A}}_\rho$
satisfies structure conditions analogous to \eqref{Eq:1:2}. Both here, and later on,
when dealing with a similar change of variable in the singular case, we drop the
suffix $\rho$ in $v_\rho$ and ${\bl{A}}_\rho$, for the sake of simplicity.

%We should now proceed with an anisotropic dyadic subdivision, in order to take into account the different sizes with respect to $x_i$ and $x_{\scriptscriptstyle N}$. 
To avoid further
technical complications, without loss of generality, we assume $L=1$.
The proof reduces to showing that
there exists a constant $\tilde\gm$ depending only on $\datap$, such
that
\begin{equation}\label{Eq:2:19}
v(x,t)\le\tilde\gm
\end{equation}
for all $(x,t)\in\{|x_i|<\frac12,\ \
0<x_{\scriptscriptstyle N}<1\}\times(-\frac{\al+\be}2c^{p-2},-\be c^{p-2}]$. In the following we denote again the
rescaled function $v$ by $u$.

Set 
\begin{equation*} 
\begin{aligned}
&K=\{|x_i|<\frac12,\ 0<x_{\scriptscriptstyle N}<1\},\\ 
&Q=K\times[-\al c^{p-2},-\be c^{p-2}].
\end{aligned}
\end{equation*}
%and operate a
% dyadic subdivision of the cube $K$ by setting $h^{\prime}=(h_1,\dots,h_{\scriptscriptstyle N-1})
% $ and considering the sub-cubes
% \begin{equation*}
% K_{h^{\prime},k}=\{(x^{\prime},x_{\scriptscriptstyle N}):\ \frac {h_i}{2^{k+1}}<x_i<\frac{h_i+1}{%
% 2^{k+1}},\ \frac1{2^{k+1}}<x_{\scriptscriptstyle N}<\frac1{2^{k}}\},
% \end{equation*}
% where $\dsty i=1,\dots,N-1,\ \ k=0,1,2,\dots,\ \ h_i=-2^{k},\dots,2^{k}-1$.
% Notice that each index $h_i$ runs independently of the others.
Suppose there exists $P_1\in Q$ such that
\begin{equation*}
u(P_1)\ge\gm^{\left\lfloor\frac{k_o}{\log_2 3/2}\right\rfloor+m},
\end{equation*}
where $\lfloor a\rfloor$ stands for the integer part of the real number $a$, $k_o\in\nn$ is
sufficiently large, $m$ will be fixed later on, and $\gm$ is the constant 
that appears in the Harnack inequality \eqref{Eq:2:7}. We claim that 
\begin{equation*}
0<x_{\scriptscriptstyle 1,N}<\left(\frac12\right)^{k_o}.
\end{equation*}
Indeed, if not, then $x_{\scriptscriptstyle 1,N}\ge\left(\frac12\right)^{k_o}$,
and by repeated application of the Harnack inequality we will show 
that this yields a contradiction. Such a procedure will also determine
the values of $\al$ and $\be$.

With respect to space variables, the worst possible case for $P_1$ is when
$x_{\scriptscriptstyle 1,N}=\left(\frac12\right)^{k_o}$, and 
$x_{\scriptscriptstyle 1,1}=\dots=x_{\scriptscriptstyle 1,N-1}=\pm\frac12$. For simplicity
let us assume $x_{\scriptscriptstyle 1,1}=\dots=x_{\scriptscriptstyle1,N-1}=\frac12$.

If we want to repeatedly apply the Harnack inequality and in this way 
getting closer and closer to $y_o=(0,\dots,0,1)$, we need to evaluate
the largest possible size of $\rho$ in \eqref{Eq:2:7} at each step.
The situation will then be the following one
\begin{equation*}
\begin{aligned}
&x_{\scriptscriptstyle 1,N}=\left(\frac12\right)^{k_o},
\quad\dist(x_1,\partial E)=\left(\frac12\right)^{k_o},\quad\rho_1=\frac12\left(\frac12\right)^{k_o},\\
&x_{\scriptscriptstyle 2,N}=\frac32\left(\frac12\right)^{k_o},
\quad\dist(x_2,\partial E)=\frac32\left(\frac12\right)^{k_o},\quad\rho_2=\frac12\,\frac32\,\left(\frac12\right)^{k_o},\\
\vdots\\
&x_{\scriptscriptstyle j,N}=\left(\frac12\right)^{k_o}\left(\frac32\right)^j,
\ \ \dist(x_j,\partial E)=\left(\frac12\right)^{k_o}\left(\frac32\right)^j,\ \ \rho_j=\frac12\left(\frac12\right)^{k_o}\left(\frac32\right)^j.\\
\end{aligned}
\end{equation*}
 We need to determine the value of $j$, at which we stop. We obviously need
 \begin{equation*}
 \left(\frac12\right)^{k_o}\left(\frac32\right)^j=1,\qquad j\approx\frac{k_o}{\log_{2}\frac32}.
 \end{equation*}
Besides getting to $y_o$, we need to have a full cube
about it, where $u$ is all bounded below, and we also need to take into account
the other coordinates, and not just $x_{\scriptscriptstyle N}$. Therefore, as $j$ needs to
belong to $\nn$, we eventually let
\begin{equation*}
j=\tilde k_o+3,
\end{equation*}
where we have set $\tilde k_o=\lfloor\frac{k_o}{\log_2 3/2}\rfloor$.
Correspondingly, by repeated application of the Harnack inequality, we conclude that
\begin{equation}\label{Eq:2:22}
\forall\,x\in K_{\frac12}(y_o)\qquad u(x,t_f)\ge\gm^{\tilde k_o+m-j}=\gm^{m-3},
\end{equation}
and now the main point becomes the evaluation of the interval
where $t_f$ can range.
When dealing with the time variable, it is easy to see that 
we have two extreme situations.

The \underline{first extreme case} is when $t_1\approx-\be c^{p-2}$ and
\begin{equation*}
u(P_1)=\gm^{\tilde k_o+m},\quad u(P_2)=\gm^{\tilde k_o+m-1},\quad\dots\quad
u(P_j)=\gm^{\tilde k_o+m-j},
\end{equation*}
that is, when the Harnack inequality gives the \underline{exact} growth of $u$.

The \underline{second extreme case} is when $t_1\approx-\al c^{p-2}$ and
\begin{equation*}
u(P_1)>>\gm^{\tilde k_o+m},\quad u(P_2)>>\gm^{\tilde k_o+m},\quad\dots\quad
u(P_j)>>\gm^{\tilde k_o+m},
\end{equation*}
that is, the function $u$ is very large and its actual decrease (if any), cannot be evaluated.

In the latter situation, we can directly assume that $t_f=-\al c^{p-2}$. Let us evaluate
what happens in the former case. By the repeated application of the Harnack inequality,
we have
\begin{equation*}
\begin{aligned}
t_f&=-\be c^{p-2}+\left(\frac{c}{\gm^{\tilde k_o+m}}\right)^{p-2}\left[\frac12\left(\frac12\right)^{k_o}\right]^p+\dots\\
&\dots+\left(\frac{c}{\gm^{\tilde k_o+m-j}}\right)^{p-2}\left[\frac12\left(\frac12\right)^{k_o}\left(\frac32\right)^j\right]^p\\
&=-\be c^{p-2}+\left(\frac{c}{\gm^{\tilde k_o+m}}\right)^{p-2}\left(\frac12\right)^{(k_o+1)p}\left(1+\gm^{p-2}\left(\frac32\right)^p+\dots\right.\\
&\left.\dots+\gm^{(p-2)j}\left(\frac32\right)^{jp}\right)\\
&=-\be c^{p-2}+\left(\frac{c}{\gm^{\tilde k_o+m}}\right)^{p-2}\left(\frac12\right)^{(k_o+1)p}\sum_{i=0}^j\left[\gm^{p-2}\left(\frac32\right)^p\right]^i\\
&=-\be c^{p-2}+\left(\frac{c}{\gm^{\tilde k_o+m}}\right)^{p-2}\left(\frac12\right)^{(k_o+1)p}
\frac{\left(\frac32\right)^{p(j+1)}\gm^{(p-2)(j+1)}-1}{\left(\frac32\right)^p\gm^{p-2}-1},
\end{aligned}
\end{equation*}
and we can conclude that
\begin{equation*}
\begin{aligned}
t_f&<-\be c^{p-2}+2\left(\frac{c}{\gm^{\tilde k_o+m}}\right)^{p-2}\left(\frac12\right)^{(k_o+1)p}
\frac{\left(\frac32\right)^{p(j+1)}\gm^{(p-2)(j+1)}}{\left(\frac32\right)^p\gm^{p-2}}\\
&<-\be c^{p-2}+2\left(\frac{c}{\gm^{\tilde k_o+m}}\right)^{p-2}\left(\frac12\right)^{(k_o+1)p}
\left(\frac32\right)^{jp}\gm^{j(p-2)}\\
&<-\be c^{p-2}+\left(\frac{c}{\gm^{\tilde k_o+m-j}}\right)^{p-2}\left(\frac94\right)^{p}=
-\be c^{p-2}+\left(\frac{c}{\gm^{m-3}}\right)^{p-2}\left(\frac94\right)^{p},
\end{aligned}
\end{equation*}
where we have taken into account the value of $j$.
Therefore, we conclude that
\begin{equation*}
t_f\in\left[-\al c^{p-2},
\left(-\be+\frac{\left(9/4\right)^p}{\gm^{(m-3)(p-2)}}\right)c^{p-2}\right].
\end{equation*}
Correspondingly, for such a $t_f$, we have not only \eqref{Eq:2:22}, but also
\begin{equation*}
\forall\,x\in K_{\frac14}(y_o)\quad u(x,t_f)\ge\gm^{m-3-l},
\end{equation*}
where $l\in\nn$ is to be fixed. The next calculations will determine $l$, and consequently $m$,
$\al$, $\be$, in order to have a contradiction and prove the claim.

We apply Lemma~\ref{Lm:2:4bis}, setting $\xi=\gm^{m-3-l}$, $a=\gm^{-1}$; the only role
played by $l$, is to provide a \emph{smaller}, initial, lower bound on $u$, and therefore 
generate a \emph{longer} cylinder $Q^+$, where the information propagates, and such that 
$(y_o,0)\in Q^+$. Since $u(y_o,0)=1$, we have a contradiction, if we end up with $u(y_o,0)>1$. 
Hence, we need to have
\begin{enumerate}[(i)]
\item $\dsty -\al<-\be$: we do not want the lower and upper bases of the cylinder $Q$ to coincide;
\item $\dsty -\be+\frac{\left(9/4\right)^p}{\gm^{(m-3)(p-2)}}\le0$: we need to be below the reference
point $(y_o,0)$ with respect to the time variable;
\item $\dsty -\al c^{p-2}+\left(\frac\dl{\gm^{m-3-l}}\right)^{p-2}\left(\frac14\right)^p\ge0$: the cylinder $Q^+$ should encompass the reference point  $(y_o,0)$;
\item $\dsty \gm^{m-4-l}>1$: this yields $u(y_o,0)>1$.
\end{enumerate}
First of all, choose $l$ such that $\dsty\gm^{(p-2)l}=\left(\frac{10^p c^{p-2}}{\dl^{p-2}}\right)$,
and then $m=l+5$: in such a way, condition (iv) is satisfied.
Notice that both $l$ and $m$ depend only on the data, but not on $k_o$.
Finally let 
\begin{equation*}
\be=\frac{\left(9/4\right)^p}{\gm^{(m-3)(p-2)}},\qquad\al= \frac{\left({10}/4\right)^p}{\gm^{(m-3)(p-2)}}.
\end{equation*}
In this way, conditions (i)--(iv) are all satisfied, and we have obtained the wanted contradiction.
Therefore, we conclude that if $u(P_1)>\gm^{\tilde k_o+m}$, then
$\dsty x_{\scriptscriptstyle 1,N}<\frac1{2^{k_o}}$. It obviously implies that 
\begin{equation}\label{Eq:2:22bis}
x_{\scriptscriptstyle 1,N}<\frac1{2^{k_o}}\qquad\text{ whenever }\qquad
u(P_1)>\gm^{m(k_o+1)}. 
\end{equation}
From here on,
we set $\gm^m=H$.
\begin{remark}
{\normalfont Instead of using Lemma~\ref{Lm:2:4bis}, the propagation of the bound below,
which generates the contradiction, can be proved by a further application of the Harnack
inequality.}
\end{remark}
%%%%%%%%%%%%%%%%%%%%%%%%%%%%%%%%%%%%%%%%%%%%%%%%
\subsubsection{\textit{End of the Proof of Theorem~\protect\ref{Thm:2:1}}}\label{S:2:4:3} 
Let
\begin{align*}
&Q_*=\left\{|x_i|<\frac14,\ 0<x_{\scriptscriptstyle N}<1,\ -\frac{\al+\be}2 c^{p-2}<t\le
-\be c^{p-2}\right\}, \\
&Q^*=\left\{|x_i|<\frac14,\ -1<x_{\scriptscriptstyle N}<1,\ -\frac{\al+\be}2 c^{p-2}<t\le
-\be c^{p-2}\right\}.
\end{align*}
Starting from $Q_*$, $Q^*$ is built by reflection; extending
$u$ to $Q^*$ as in Lemma~\ref{Lm:2:4}, $u$ is
still a (signed) solution to \eqref{Eq:1:1}--\eqref{Eq:1:2}. 

\noi Now let $\dsty %
P_1=(x_1,t_1)=(x_{\scriptscriptstyle 1,1},x_{\scriptscriptstyle 1,2},\dots,x_{\scriptscriptstyle 1,N},t_1) \in Q_*$ be such that
\begin{equation}  \label{Eq:2:23}
u(P_1)\ge H^{k_o+1};
\end{equation}
by \eqref{Eq:2:22bis}, it must be $0<x_{\scriptscriptstyle 1,N}<\frac1{2^{k_o}}$. Set
\begin{equation*}
Q(P_1)=\left\{|x-x_1|<2^{-k_o}\eps^{-s},\ \ t_1-(2^{-k_o}\eps%
^{-s})^p<t<t_1\right\},
\end{equation*}
where $\eps$ is the quantity claimed by Lemma~\ref{Lm:2:3} and $s\in\nn$ is
to be fixed. Let
\begin{equation*}
\om_o^{(1)}=\osc_{Q(P_1)}u.
\end{equation*}
Without loss of generality, we may assume $\eps<\frac12$. Thanks to \eqref{Eq:2:23} and the
construction of $u$ by odd reflection, we have
\begin{equation}  \label{Eq:2:24}
\om_o^{(1)}\ge 2H^{k_o+1}.
\end{equation}
Moreover, if $k_o$ is large enough, we have $Q(P_1)\subset Q^*$. Set $\sig%
_o^{(1)}=2^{-k_o}\eps^{-s}$ and consider
\begin{equation*}
Q_o^{(1)}=K_{\sig_o^{(1)}}(x_1)\times(t_1-\theta_o^{(1)}(\sig%
_o^{(1)})^p,t_1), \qquad\text{where}\ \ \ \theta_o^{(1)}=\left(\frac c{\om%
_o^{(1)}}\right)^{p-2}.
\end{equation*}
It is apparent that $Q_o^{(1)}\subset Q(P_1)\subset Q^*$. Notice that we do
not need to assume \eqref{Eq:2:13} here, since, by construction, the
cylinders are all correctly nested into one another. By Lemma~\ref{Lm:2:3},
we can build a sequence
\begin{equation*}
\om_n^{(1)}=\dl\om_{n-1}^{(1)},\quad \theta_n^{(1)}=\left(\frac c{\om_n^{(1)}}\right)^{p-2},\quad
\sig_n^{(1)}=\eps\sig_{n-1}^{(1)},\quad Q_n^{(1)}=Q_{\sig_n^{(1)}}(\theta_n^{(1)}),
\end{equation*}
for all non-negative integers $n$. Such a sequence satisfies
\begin{equation*}
Q_{n+1}^{(1)}\subset Q_n^{(1)},\qquad \osc_{Q_n^{(1)}}u\le\om_n^{(1)}.
\end{equation*}
By iteration
\begin{equation*}
\osc_{Q_n^{(1)}}u\le\dl^n\om_o^{(1)}=\dl^n\osc_{Q(P_1)}u,
\quad\Rightarrow\quad\osc_{Q(P_1)}u\ge\frac1{\dl^n}\osc_{Q_n^{(1)}}u.
\end{equation*}
If we now choose $n=s$, and $s$ such that $\dl^{-s}>H^{10}$, by the choice
of $\sig_o^{(1)}$ we conclude that
\begin{equation*}
\om_o^{(1)}\ge2H^{k_o+11}
\end{equation*}
and this obviously improves the previous lower bound given by \eqref{Eq:2:24}. 
As $u$ has been built by odd reflection, we conclude there must exist $%
P(x_2,t_2)=(x_{\scriptscriptstyle 2,1},x_{\scriptscriptstyle 2,2},\dots,x_{\scriptscriptstyle 2,N},t_2)\in Q(P_1)$ such that
\begin{equation*}
u(P_2)\ge H^{k_o+11}.
\end{equation*}
As before, by \eqref{Eq:2:22bis}, $0<x_{\scriptscriptstyle 2,N}<\frac1{2^{k_o+10}}$, and also $%
t_1-(2^{-k_o}\eps^{-s})^p<t_2<t_1$. Set
\begin{equation*}
Q(P_2)=\left\{|x-x_2|<2^{-k_o-10}\eps^{-s},\ \ t_2-(2^{-k_o-10}\eps%
^{-s})^p<t<t_2\right\}.
\end{equation*}
Once more, provided $k_o$ is large enough, we can assume that $Q(P_2)\subset
Q^*$. Arguing as before, we conclude there exists $P_3(x_3,t_3)\in
Q(P_2)$ such that
\begin{equation*}
u(P_3)\ge H^{k_o+21}.
\end{equation*}
Consequently $0<x_{\scriptscriptstyle 3,N}<2^{-k_o-20}$. By induction, we get 
$\{P_q(x_q,t_q)\}$, such that
\begin{equation}  \label{Eq:2:25}
u(P_q)\ge H^{k_o+1+10(q-1)},
\end{equation}
and
\begin{equation*}
0<x_{\scriptscriptstyle q,N}<2^{-k_o-10(q-1)}.
\end{equation*}
Notice that $k_o$ depends on $\al$, $\be$, $c$, and $s$, and therefore,
due the definition of these quantities, on the data $p$, $N$, $C_o$, $C_1$.
Now choose $k_o$ so large as to have
\begin{equation*}
t_1-\sum_{q=1}^{\infty}(2^{-k_o-10(q-1)}\eps^{-s})^p>-\al c^{p-2},
\end{equation*}
and $\forall\,i=1,\dots,N-1$
\begin{equation*}
x_{\scriptscriptstyle 1,i}-\sum_{q=1}^\infty(2^{-k_o-10(q-1)}\eps^{-s})>-\frac12,\quad
x_{\scriptscriptstyle 1,i}+\sum_{q=1}^\infty(2^{-k_o-10(q-1)}\eps^{-s})<\frac12.
\end{equation*}
It is the need to satisfy these requirements that forces $|x_i|<\frac14$ in the definition of
$Q_*$ and $Q^*$. Indeed the value of $k_o$ determines the value of $\tilde\gm$
in \eqref{Eq:2:19}, and therefore, we must be able to choose $k_o$, independently
of the solution or any other geometrical conditions. If $|x_i|$ were larger than $\frac14$,
then $k_o$ would also depend on the distance in space of $P_1$
to the boundary of the cube of edge $\frac12$.

Once the previous conditions are satisfied,
the sequence $\{P_q\}$ is contained in a fixed cylinder of $Q_*$.
Together with \eqref{Eq:2:25}, this leads to a contradiction, since the sequence
approaches the boundary, and the corresponding values of $u$ grow arbitrarily large, 
whereas $u$ is assumed to vanish continuously at the boundary.
\hfill\bbox  %
\vskip.2truecm

\subsubsection{\textit{Proof of Corollary~\protect\ref{Cor:2:1}}}

Let
\begin{equation*}
\mathcal{K}=\overline{\Psi^-_{\frac\rho2}\pto}\cap(\partial E\times(0,T)),\qquad \dsty M=\sup_{\Psi^-_{\frac\rho2}\pto}u.
\end{equation*}
By Theorem~1.2 of Chapter III of \cite{dibe-sv}, $\forall(x,t)\in{%
\Psi^-_{\frac\rho2}\pto}$ we have
\begin{equation*}
u(x,t)\le \gm M\left(\frac{\inf_{(y,s)\in\mathcal{K}}\left(|x-y|+M^{\frac{p-2%
}p}|t-s|^{\frac1p}\right)}\rho \right)^{\mu}.
\end{equation*}
Since we are dealing with 
a Lipschitz cylinder, in the infimum above we can take $t\equiv s$, and we reduce to
\begin{equation*}
u(x,t)\le \gm M\left(\frac{\dist(x,\partial E)}\rho \right)^{\mu}.
\end{equation*}
By \eqref{Eq:2:1} and
possibly a further application of the Harnack inequality, we conclude. \hfill%
\bbox

\subsubsection{\textit{Proof of Theorem~\protect\ref{Thm:2:2}}}

Since the boundary is of class $C^2$, 
the portion of the boundary $\partial E\cap K_{2\rho}(x_o)$ can now be represented by
$\dsty x_{\scriptscriptstyle N}=\Phi(x^{\prime})$, where $\Phi$ is a function of class $C^2$
satisfying
\begin{equation*}
\Phi(x_o')=0,\qquad D\Phi(x_o')=0,\qquad\| D^\prime\Phi\|_\infty\le
M_2,\qquad\| (D^\prime)^2\Phi\|_\infty\le M_2,
\end{equation*}
and $M_2\in(0,1)$ is a proper parameter, provided $\rho$ is small enough.

\noi If we now introduce the new variables
\begin{equation*}
y_i=x_i,\quad i=1,\dots,N-1,\qquad y_{\scriptscriptstyle N}=x_{\scriptscriptstyle N}-\Phi(x^{\prime}),
\end{equation*}
the portion of the boundary $\partial E\cap K_{2\rho}(x_o)$ coincides with the
portion of the hyperplane $y_{\scriptscriptstyle N}=0$ within $K_{2\rho}(y_o)$. We orient $y_{\scriptscriptstyle N}$ so
that $\dsty E\cap K_{2\rho}(y_o)\subset\{y_{\scriptscriptstyle N}>0\}$. Denoting again with $x$ the
transformed variable $y$, we proceed as in the proof of Theorem~\ref{Thm:2:1}. 
Consequently (1.1)$_o$ is rewritten as
\begin{equation}  \label{Eq:2:26}
u_t-\dvg{\bl A}(x,Du)=0,
\end{equation}
%where, 
%
%We have
%\begin{equation*}
%{\bl A}(x,\etavect)=[\langle\etavect,\overline{\mathbb{B}}\etavect%
%\rangle]^{(p-2)/2}\overline{\mathbb{B}}\etavect =\sum_{k=1}^Na_{ik}(x,%
%\etavect)\eta_k,
%\end{equation*}
{where, just as in the proof of Theorem~\ref{Thm:2:1},
\begin{equation}  \label{Eq:2:27}
\left\{
\begin{array}{l}
\bl{A}(x,\etavect)\cdot \etavect\ge C_o|\etavect|^p \\
|\bl{A}(x,\etavect)|\le C_1|\etavect|^{p-1}%
\end{array}%
\right .\quad \text{ a.e. }\> (x,t)\in E_T,
\end{equation}
$C_o$, $C_1$ are positive constants depending only on $N$, $p$, $M_2$,
%\[
%\mathbb{B}:=\left[
%\begin{matrix}
%\mathbb{I}_{\scriptscriptstyle N-1} &- D^\prime_x\Phi\cr & \cr 0 & 1%
%\end{matrix}
%\right],
%\]
%and
%\begin{equation*}
%\begin{aligned}
%&[a_{ik}(x,\etavect)]=\Big[|\etavect|^2+\eta_{\scriptscriptstyle N}^2| D^\prime\Phi(x)|%^2-2\eta_{\scriptscriptstyle N}%
%\langle\etavect^{\prime},  D^\prime\Phi(x)\rangle\Big]^{\frac{p-2}2}\\
%&\times \left[
%\begin{matrix}
%\mathbb{I}_{\scriptscriptstyle N-1} & - D^\prime\Phi(x)\cr - D^\prime\Phi(x) & 1+| D^\prime\Phi(x)|%^2\cr%
%\end{matrix}
%\right].
%\end{aligned}
%\end{equation*}
and
\begin{equation*}
\etavect^{\prime}=(\eta_1,\eta_2,\dots,\eta_{\scriptscriptstyle N-1}),\qquad |\etavect%
|^2=\sum_{k=1}^N|\eta_k|^2.
\end{equation*}
Our solution to \eqref{Eq:2:26} vanishes continuously on $%
K_{2\rho}(x_o)\cap\{x_{\scriptscriptstyle N}=0\}$.
Moreover, by lengthy but quite straightforward calculations,
\begin{align} 
&\frac{\partial A_i}{\partial\eta_k}\xi_i\xi_k\ge C_2|\etavect|^{p-2}|\xi|^2, \label{Eq:2:28}\\
&\left|\frac{\partial A_i}{\partial\eta_k}\right|\le C_3|\etavect%
|^{p-2},\qquad \left|\frac{\partial A_i}{\partial x_k}\right|\le C_4|\etavect%
|^{p-1},\label{Eq:2:29}
\end{align}
where $C_2$, $C_3$, $C_4$ depend only on $N$, $p$, and $M_2$.} Finally, due to its
definition, ${\bl A}$ satisfies the homogeneity condition
\begin{equation}  \label{Eq:2:30}
{\bl A}(x,\etavect)=|\etavect|^{p-1}{\bl A}(x,\frac{\etavect}{|\etavect|}).
\end{equation}
Therefore, proving Theorem~\ref{Thm:2:2} reduces to proving
\begin{lemma}%\label{Lm:2:6}
Let $u$ be a non-negative, weak solution to
\[
u_t-\dvg{\bl A}(x,Du)=0\qquad\text{in}\ E_T
\]
for $p>2$, where ${\bl A}$ satisfies the structure conditions \eqref{Eq:2:27}, \eqref{Eq:2:28},
\eqref{Eq:2:29}, and \eqref{Eq:2:30}.
Take $\pto\in S_T$, $\rho\in(0,r_o)$, let $P_\rho=(x_o',2M_2\rho,t_o)$,
and assume that $u\prho>0$, $\partial E$ is flat with respect to $x_{\scriptscriptstyle N}$, and 
$(t_o-\theta(4\rho)^p,t_o+\theta(4\rho)^p]\subset(0,T]$,
where $\theta=\left[\frac c{u\prho}\right]^{p-2}$.

\noi Suppose that $u$ vanishes continuously 
on $(\partial E\cap K_{2\rho}(x_o))\times(t_o-\theta(4\rho)^p,t_o+\theta(4\rho)^p]$.
Then there exists a constant $\gm>0$, depending only on $N$, $p$,
and  $C_i$, $i=0,\dots,4$, such that
\begin{equation*}%%\label{Eq:2:31}
0\le u(x',x_{\scriptscriptstyle N},t)\le\gm\left(\frac{x_{\scriptscriptstyle N}}\rho\right) u\prho,
\end{equation*}
for all
\begin{equation*}
(x,t)\in\{|x_i-x_{o,i}|<\frac{\rho}2,\ 0<x_{\scriptscriptstyle N}< 2M_2\rho\}\times(t_o-\frac{\al+3\be}4\theta\rho^p,t_o-\be\theta\rho^p].
\end{equation*}
\end{lemma}
\vskip.2truecm
\noi{\it Proof} - We use the same argument of \cite[Theorem~4.1]{DKV}. By possible,
suitable rescaling and translation, assume $x_o=0$, $\rho=1$, and 
let ${\cal M}=\sup_{\Psi^-_1}u$.
By \eqref{Eq:2:18}, ${\cal M}\le\tilde\gm\, u(P_1)$. Let $y=(0,\dots,0,-1)$,
$t_1=t_o-\be\theta$, and consider the function
\begin{equation*}
\eta_k(x,t)=\exp[-k(|x-y|-1)]\exp[u(P_1)^{p-2}(t-t_1)],
\end{equation*}
and the set
\begin{equation*}
{\cal N}_k=\{(x,t):\,x_{\scriptscriptstyle N}>0,\ 1<|x-y|<1+\frac1k,\ t_1-(1+\dl)\be\theta<t<t_1\},
\end{equation*}
for some small enough, positive parameter $\dl$.
We assume $k$ is so large, that ${\cal N}_k\subset\Psi^-_1$.
If we choose $\dsty C=\tilde\gm\max\{(1-e^{-1})^{-1};(1-e^{-(1+\dl)\be c^{p-2}})^{-1}\}$, 
where $\tilde\gm$ is the constant of \eqref{Eq:2:1}, it is easy to verify that 
$\dsty \Theta_k(x,t)=C u(P_1)(1-\eta_k(x,t))$ satisfies $u\le\Theta_k$ on the parabolic 
boundary of ${\cal N}_k$.
Provided we choose $k$ as the largest positive root of
\begin{equation*}
C_o(p-1)k^p - b k^{p-1}- C^{2-p}=0,
\end{equation*}
where $b$ is a positive quantity that depends only on $M_2$, and $C_o$ is the constant 
of the first of \eqref{Eq:2:27}, relying on \eqref{Eq:2:27}--\eqref{Eq:2:30},
it is a matter of straightforward calculations, to verify that $\Theta_k$ is a super-solution
to \eqref{Eq:2:26} in ${\cal N}_k$.

\noi By the comparison principle, $u\le\Theta_k$ in ${\cal N}_k$. In particular,
$\forall\,0<x_{\scriptscriptstyle N}<\frac1k$,
\begin{equation*}
\begin{aligned}
u(0,\dots,0,x_{\scriptscriptstyle N},t_1)&\le\Theta_k(0,\dots,0,x_{\scriptscriptstyle N},t_1)\\
&=C u(P_1)(1-e^{-kx_{\scriptscriptstyle N}})\\
&\le\tilde\gm\, x_{\scriptscriptstyle N}\,u(P_1).
\end{aligned}
\end{equation*}
On the other hand, if $x_{\scriptscriptstyle N}\ge\frac1k$, then
\begin{equation*}
u(0,\dots,0,x_{\scriptscriptstyle N},t_1)\le \mathcal M\le\tilde\gm\,k\,x_{\scriptscriptstyle N}\, u(P_1).
\end{equation*}
The same argument can be repeated using any $t_1\in(t_o-\frac{\al+3\be}4\theta,t_o-\be\theta]$,
and switching back to the original coordinates, we conclude.
\hfill\bbox
\begin{remark}
{\normalfont Instead of relying on the Comparison Principle, Theorem~\ref{Thm:2:2}
could be proved using the $L^\infty$-estimates for the gradient $Du$, as shown, for
example, in \cite[Chapter~VIII]{dibe-sv} (see also \cite{KM2,KM3}). However, similar estimates
are not known for solutions to equation \ref{Eq:4:1:o}, and one does not expect them to
hold true for solutions to \eqref{Eq:4:1}--\eqref{Eq:4:2}: therefore, we
opted for an approach that works in both cases.}
\end{remark}
\section{The Singular Super-critical Case $\frac{2N}{N+1}<p<2$}\label{S:3}

\subsection{A Weak Carleson Estimate}
We consider our first result in the singular super-critical case $\frac{2N}{N+1}<p<2$.  
Let $E_T$, $u$, $\pto$, $\rho$, $x_\rho$, $P_\rho$ be as in Theorem~\ref{Thm:2:1} and
set 
\begin{equation*}%%\label{Eq:3:1}
I(t_o,\rho,h)=(t_o- h \rho^p,t_o+ h \rho^p).
\end{equation*}
Moreover, let $u$ be a weak solution to \eqref{Eq:1:1}--\eqref{Eq:1:2}
such that 
\begin{equation}\label{Eq:3:1bis}
0<u\le M\qquad\text{ in }\ \ E_T,
\end{equation} 
and assume that
\begin{equation}\label{Eq:3:1ter}
I(t_o,9\rho,M^{2-p})\subset(0,T].
\end{equation}
Then we define
\begin{equation*}
\begin{aligned}
&\widetilde\Psi_{\rho}\pto=E_T\cap\left\{(x,t):\ |x_i-x_{o,i}|<2\rho,\ |x_{\scriptscriptstyle N}|<4L\rho,
t\in I(t_o,9\rho,\etarho^{2-p})\right\},\\
&\bar\Psi_{\rho}\pto=E_T\cap\left\{(x,t):\ |x_i-x_{o,i}|<\frac\rho4,\ |x_{\scriptscriptstyle N}|<2L\rho,
t\in I(t_o,\rho,\etarho^{2-p})\right\},
\end{aligned}
\end{equation*}
where $\etarho$ is the first root of the equation
\begin{equation}  \label{Eq:3:1quater}
\max_{\widetilde\Psi_{\rho}\pto}u=\etarho.
\end{equation}
Notice that both the functions
%\begin{equation*}
$\dsty y_1(\etarho)=\max_{\widetilde\Psi_{\rho}\pto}u$, 
$\dsty y_2(\etarho)=\etarho$
%\end{equation*}
are monotone increasing. Moreover
\begin{equation*}
\left\{
\begin{array}{l}
y_1(0)\ge u(P_\rho)>0, \\
y_2(0)=0,%
\end{array}%
\right .\quad \text{ and }\quad \left\{
\begin{array}{l}
y_1(M)\le M, \\
y_2(M)=M.%
\end{array}%
\right .
\end{equation*}
Therefore, it is immediate to conclude that at least one root of \eqref{Eq:3:1quater}
actually exists. Moreover, by \eqref{Eq:3:1ter} $\widetilde\Psi_{\rho}\pto\subset E_T$. 

\noi As already mentioned in the Introduction, we  can only provide 
a weak form of the Carleson estimate, expressed by the following theorem.

\begin{theorem}\label{Thm:3:1}
\emph{(Carleson-type Estimate, weak form, $\frac{2N}{N+1}<p<2$).} Let $u$ 
be a weak solution to \eqref{Eq:1:1}--\eqref{Eq:1:2}, that satisfies 
\eqref{Eq:3:1bis}. Assume that \eqref{Eq:3:1ter} holds true and $u$ vanishes continuously on
$$\partial E\cap \{|x_i-x_{o,i}|<2\rho,|x_{\scriptscriptstyle N}|<8L\rho\}\times I(t_o,9\rho,M^{2-p}).$$
Then there exist constants $\gm>0$ and $\al\in(0,1)$, depending only on
 $\datap$ and $L$, such that
\begin{equation*}
u(x,t)\le\gm\left(\frac{\dist(x,\partial E)}{\rho}\right)^{\al}
\times\sup_{\tau\in I(t_o,\rho,2\etarho^{2-p})}u(x_\rho,\tau),
\end{equation*}
for every $(x,t)\in\bar\Psi_{\rho}\pto$.
\end{theorem}
If we let
\begin{equation*}
\Psi_{\rho,M}\pto=E_T\cap\left\{(x,t):\ |x_i-x_{o,i}|<\frac\rho4,\ |x_{\scriptscriptstyle N}|<2L\rho,
t\in I(t_o,\rho,M^{2-p})\right\},
\end{equation*}
we have a second statement.
\begin{corollary}\label{Cor:3:01}
Under the same assumptions of Theorem~\ref{Thm:3:1}, we have
\begin{equation*}
u(x,t)\le\gm\left(\frac{\dist(x,\partial E)}{\rho}\right)^{\al}
\times\sup_{\tau\in I(t_o,\rho,2M^{2-p})}u(x_\rho,\tau),
\end{equation*}
for every $(x,t)\in\Psi_{\rho,M}\pto$.
\end{corollary}

The quantity $\etarho$ is known only qualitatively through \eqref{Eq:3:1quater},
whereas $M$ is a datum. Therefore, Corollary~\ref{Cor:3:01} 
can be viewed as a quantitative version of a purely qualitative statement. On
the other hand, since $\etarho$ could be attained in $P_\rho$, Theorem~\ref{Thm:3:1} gives the
sharpest possible statement, and is genuinely intrinsic.

Moreover, with respect to Theorem~\ref{Thm:2:1} and Corollary~\ref{Cor:2:1}, 
Theorem~\ref{Thm:3:1} combines two distinct statements in a single one
(mainly for simplicity), and presents two fundamental differences:
when $p>2$, the value of $u$ at a point above,
controls the values of $u$ below, whereas when $\frac{2N}{N+1}<p<2$, 
the {\bf maximum} of $u$ over a proper time interval centered at $t_o$
controls the values of $u$ {\bf both above and below} the time level $t_o$. 
These are consequences of the different statements of the Harnack inequality 
in the two cases. In fact, the following theorem is proved in \cite{DBGV-sing07} 
(see also \cite{DBGV-mono} for a thorough presentation).

For fixed $(x_o,t_o)\in E_T$ and $\rho>0$, set
%\begin{equation*}
$\dsty {\mathcal M}=\sup_{K_{\rho}(x_o)}u(x,t_o)$,
%\end{equation*}
and require that
\begin{equation} \label{Eq:3:2}
K_{8\rho}(x_o)\times I(t_o,8\rho,{\mathcal M}^{2-p}) \subset
E_T.
\end{equation}
\begin{theorem}\label{Thm:3:2}
\emph{(Harnack Inequality)} Let $u$ be a
non-negative, weak solution
to \eqref{Eq:1:1}--\eqref{Eq:1:2},
in $E_T$ for $p\in(\frac{2N}{N+1},2)$.

There exist constants $\overline{\eps}\in(0,1)$ and $\overline{\gm}>1$,
depending only on $\datap$, such
that for all intrinsic cylinders $\pto+Q_{8\rho}^{\pm}(\theta)$
for which \eqref{Eq:3:2} holds,
\begin{equation}\label{Eq:3:3}
\overline{\gm}^{-1}\sup_{K_{\rho}(x_o)} u(\cdot,\sig)\le u\pto\le
\overline{\gm} \inf_{K_{\rho}(x_o)} u(\cdot,\tau)
\end{equation}
for any pair of time levels $\sig,\tau$ in the range
\begin{equation}\label{Eq:3:4}
t_o-\overline{\eps}\, u(x_o,t_o)^{2-p}\rho^p\le\,\sig,\tau\,
\le t_o+\overline{\eps}\, u(x_o,t_o)^{2-p}\rho^p.
\end{equation}
The constants $\overline{\eps}$ and $\overline{\gm}^{-1}$ tend to zero as either
$p\to2$ or as $p\to\frac{2N}{N+1}$.
\end{theorem}
%%%%%%%%%%%%%%%%%%%%%%%%%%%%%%%%%%%%%%%%%%%%%%%%
\begin{remark}
{\normalfont  With respect to the degenerate case, we now have $c=1$ for the size of the intrinsic cylinders. The upper bound ${\mathcal M}$ has only the qualitative role to insure that $\pto+Q^{\pm}_{8\rho}({\mathcal M})$ are contained
within the domain of definition of $u$.}
\end{remark}
%%%%%%%%%%%%%%%%%%%%%%%%%%%%%%%%%%%%%%%%%%%%%%%%
\subsection{A Counterexample}

Can we improve the result of Theorem~\ref{Thm:3:1}, namely can we substitute
the supremum of $u$ on $I(t_o,\rho,2\etarho^{2-p})$ with the pointwise value $u\prho$? 
This would certainly be possible, if there existed a constant $\gm$, dependent only on the data
$\datap$, such that
\begin{equation*}
\forall\,t\in I(t_o,\rho,2\etarho^{2-p})\qquad u(x_\rho,t)\le\gm\, u\prho.
\end{equation*}
Under a geometrical point of view, this amounts to building
a Harnack chain connecting $(x_\rho,t)$ and $P_\rho$, for all $t\in I(t_o,\rho,2\etarho^{2-p})$. 
In general, 
without further assumptions on $u$, 
this is not possible, as the following counterexample shows.

\noi Let $u$ be the unique non-negative solution to
\begin{equation*}
\left\{
\begin{array}{l}
u\in C(\rr_+;L^2(E))\cap L^p(\rr_+;W^{1,p}_o(E))\\
u_t-\dvg(|Du|^{p-2}Du)=0\quad\text{ in }\,\, E_T\\
u(\cdot,0)=u_o\in C^o(\overline E),\\
\end{array}%
\right .
\end{equation*}
with $u_o>0$ in $E$, and $u_o\big|_{\partial E}=0$.

\noi By Proposition~2.1, Chapter~VII of \cite{dibe-sv}, there exists a finite time $T_*$, depending
only on $N$, $p$, $u_o$, such that $u(\cdot,t)\equiv0$ for all $t\ge T_*$. 
By the results of \cite[Chapter~IV]{dibe-sv}, $u\in C^o(\overline{E\times(0,T_*)})$. Suppose
now that at time $t=T_*+1$, we modify the boundary value and for any $t>T_*+1$ we let
$\dsty u(\cdot,t)\big|_{\partial E}=g(\cdot,t)$, where $g$ is continuous and strictly positive.
It is immediate to verify that $u$ becomes strictly positive for any $t>T_*+1$. Therefore, the
positivity set for $u$ is not a connected set, $\dsty u(x,t)\equiv0$ for all $\dsty\forall\,(x,t)\in
\overline{E\times(T_*,T_*+1)}$, and if $(x_\rho,t)$ and $P_\rho$ lie on opposite sides of the
vanishing layer for $u$, by the intrinsic nature of Theorem~\ref{Thm:3:2}, 
there is no way to connect them with a Harnack chain.

The previous counterexample allows $u$ to vanish identically for $t$ in a proper interval, 
but by suitably modifying the boundary values, it is clear that we can have $u$ strictly positive, 
and as close to zero as we want. Therefore, the impossibility of connecting two arbitrary points by a Harnack chain, does not depend on the vanishing of $u$, but it is a general property of solutions
to \eqref{Eq:1:1}--\eqref{Eq:1:2}, whenever $E\not\equiv\rn$. Moreover, by properly adjusting the boundary value, one can even create an arbitrary number of oscillations for $u$ between 
positivity and null regions.

We considered solutions to the $p$--Laplacian just for the sake of simplicity, 
but everything continues to hold, if we consider the same boundary value problem
for \eqref{Eq:1:1}--\eqref{Eq:1:2}.

Notice that if we deal with weak solutions to \eqref{Eq:1:1}--\eqref{Eq:1:2}
in $\rn\times(0,T]$, then we do not have boundary values any more, the situation
previously discussed cannot occur, and therefore
any two points $(x,t)$ and $\pto$ can always be connected by a Harnack chain, 
provided both $u(x,t)$ and $u\pto$ are strictly positive, and
$0<t-t_o<\frac{\eps}{8^p}t_o$, as discussed in \cite[Chapter~7, Proposition~4.1]{DBGV-mono}. 
The sub potential lower bound discussed there
is then a property of weak solutions given in the whole $\rn\times(0,T)$.

The Harnack inequality given in Theorem~\ref{Thm:3:2} is time-insensitive,
and its constants are not stable as $p\to2$.
A different statement, analogous to the one given in Theorem~\ref{Thm:2:3}, could be given,
and in such a case the constants would be stable (see \cite[Chapter~6]{DBGV-mono} for a thorough discussion of the two possible forms). However, the eventual result is the same, 
and independently of the kind of Harnack inequality one considers, two points
$(x,t)$ and $\pto$ of positivity for $u$, cannot be connected by a Harnack chain. 

Notice that we have a sort of dual situation: when $1<p<2$ the support of $u$ can be disconnected in time, when $p>2$, as we discussed in \S~\ref{SS:2:2}, the support can be disconnected in space.

Strictly speaking, the previous counterexample only shows that we cannot replace the line with a point, but per se it does not rule out the possibility for a strong form of Carleson's estimate to hold true
all the same. 
However, if one tries to adapt to the singular super-critical context the standard proof 
based on the Harnack inequality and the boundary H\"older continuity (as we did, for example, 
in the degenerate context), then one quickly realizes that, in order to have the
cylinders of Lemma~\ref{Lm:2:3} inside the reference cylinder, one needs to 
know in advance the oscillation of $u$: this suggests that only a control in terms of the 
supremum taken in a proper set can be feasible.
%%%%%%%%%%%%%%%%%%%%%%%%%%%%%%%%%%%%%%%%%%%%%%%
\subsection{A Strong Carleson Estimate}\label{SS:strongCarleson}

With respect to the statement of Theorem~\ref{Thm:3:1}, 
a stronger form is indeed possible, provided we allow the parameter $\gm$ 
to depend not only on the data, but also on the oscillation of $u$. 

\noi Let $E_T$, $u$, $\pto$, $\rho$, $P_\rho$ be as in Theorem~\ref{Thm:2:1}, 
and for $k=0,1,2,\dots$ set
\begin{equation*}
\begin{aligned}
&\rho_k=\left(\frac78\right)^k\rho,\qquad\sig_k=\frac{\rho_k}{\gm^{k\frac{2-p}p}},\\
&x_{\rho_k}=(x_o',2L\rho_k),\qquad P_{\rho_k}=(x_o',2L\rho_k,t_o),\\
&\Psi_{\rho_k,M}\pto\\
&\quad=E_T\cap\left\{(x,t):|x_i-x_{o,i}|<{\frac{\rho_k}4},\
|x_{\scriptscriptstyle N}|<2L\rho_k, \ t\in I(t_o,\sig_k,M^{2-p})\right\},\\
&m_o=\inf_{\tau\in I(t_o,\rho,2M^{2-p})}u(x_\rho,\tau),\quad 
M_o=\sup_{\tau\in I(t_o,\rho,2M^{2-p})}u(x_\rho,\tau).
\end{aligned}
\end{equation*}

\begin{corollary}\label{Cor:3:1}
\emph{(Carleson-type Estimate, strong form, $\frac{2N}{N+1}<p<2$)}. Let $u$ be a
weak solution to \eqref{Eq:1:1}--\eqref{Eq:1:2} 
such that $0<u\le M$ in $E_T$. Assume that
$I(t_o,9\rho,M^{2-p})\subset(0,T]$ and that $u$ vanishes continuously on
$$\partial E\cap \{|x_i-x_{o,i}|<2\rho,|x_{\scriptscriptstyle N}|<8L\rho\}\times I(t_o,9\rho,M^{2-p}).$$
Then there exists a constant $\gm$, depending only on
 $\datap$, $L$, and $\frac M{m_o}$, such that
\begin{equation}\label{Eq:3:7}
u(x,t)\le\gm\,u(P_{\rho_k}),
\end{equation}
for every $(x,t)\in\Psi_{\rho_k,M}\pto$, for all $k=0,1,2,\dots$.
\end{corollary}

\begin{remark}
{\normalfont The strong form of the Carleson-type estimate is derived from Corollary~\ref{Cor:3:01}.
An analogous statement can be derived from Theorem~\ref{Thm:3:1}.}
\end{remark}

\begin{remark}
{\normalfont Estimate \eqref{Eq:3:7} has the same structure
as the backward Harnack inequality for caloric
functions that vanish just on a disk at the boundary (see \cite[Theorem~13.7, page 234]{CS}).
This is not surprising, because \eqref{Eq:3:7} is indeed a backward Harnack 
inequality, due to the specific nature of the Harnack inequality for the singular case.
However, it is worth mentioning that things are not completely equivalent;
indeed, the constants we have in the time-insensitive Harnack inequality 
\eqref{Eq:3:3}--\eqref{Eq:3:4} are not stable (and cannot be stabilized), and therefore, 
the result for caloric functions cannot be recovered from the singular case, 
by simply letting $p\to2$ (as it is instead the case for many other results).}
\end{remark}
\vskip.2truecm
\noi Another striking difference with respect to the degenerate case, 
appears when we consider $C^{1,1}$ cylinders and (mainly for simplicity) 
the prototype equation \ref{Eq:1:1:o}. In this case, indeed, weak solutions
vanishing on the lateral part enjoy a linear behavior at the boundary with 
implications expressed in the following result. Note that the role of $L$
in the definition of $\Psi_{\rho,M}$ is now played by $M_{1,1}$.
\begin{theorem}\label{Thm:3:3}
Let $\frac{2N}{N+1}<p<2$. Assume $E_T$ is a $C^{1,1}$ cylinder,
and $\pto$, $\rho$, $P_\rho$ as in Theorem~\ref{Thm:2:1}. 
Let $u,\,v$ be two weak solutions to \ref{Eq:1:1:o} in $E_T$, 
satisfying the hypotheses of Theorem~\ref{Thm:3:1}, 
$0<u,v\le M$ in $E_T$.
Then there exist positive  constants $\bar s$, $\gm$, $\be$, $0<\be\le1$,
depending only on $N$,
$p$, and $M_{1,1}$, and $\rho_o$, $c_o>0$, depending also on the oscillation of $u$,
such that the following properties hold.
\begin{description}
\item[(a)] \emph{Hopf Principle:}
\begin{equation}\label{Eq:3:8}
|Du|\ge c_o\qquad\text{in}\qquad  {\Psi_{\rho_o,M}(x_o,t_o)}.
\end{equation}
\item[(b)] \emph{Boundary Harnack Inequality:} 
\begin{equation}\label{Eq:3:8bis}
\gm^{-1}\,\frac{\dsty\inf_{\tau\in I(t_o,\rho,2M^{2-p})}u(x_\rho,\tau)}{\dsty\sup_{\tau\in I(t_o,\rho,2M^{2-p})}v(x_\rho,\tau)}\le \frac{u(x,t)}{v(x,t)}\le\gm\,\frac {\dsty\sup_{\tau\in I(t_o,\rho,2M^{2-p})}u(x_\rho,\tau)}{\dsty\inf_{\tau\in I(t_o,\rho,2M^{2-p})}v(x_\rho,\tau)},
\end{equation}
for all $(x,t)\in\{x\in K_{\bar s\frac\rho4}(x_o)\cap E:\,\dist(x,\partial E)<\bar s\frac\rho8\}\times I(t_o,\rho,\frac12M^{2-p})$,
with $\rho<\rho_o$.
 \item[(c)]  The quotient $u/v$ is H\"older continuous with exponent $\be$ in 
 ${\Psi_{\frac{\rho_o}2,M}(x_o,t_o)}$
\end{description}
\end{theorem}
%%%%%%%%%%%%%%%%%%%%%%%%%%%%%%%%%%%%%%%%%%%%%%%%
\begin{remark}\label{Rmk:3:4}
{\normalfont Since
\[
\frac {\dsty\sup_{\tau\in I(t_o,\rho,2M^{2-p})}u(x_\rho,\tau)}{\dsty\inf_{\tau\in I(t_o,\rho,2M^{2-p})}v(x_\rho,\tau)}\le\frac{M_{o,u} u(P_\rho)}{m_{o,u}}\frac{M_{o,v}}{m_{o,v}v(P_\rho)},
\]
\[
\frac{\dsty\inf_{\tau\in I(t_o,\rho,2M^{2-p})}u(x_\rho,\tau)}{\dsty\sup_{\tau\in I(t_o,\rho,2M^{2-p})}v(x_\rho,\tau)}\le\frac{m_{o,u} u(P_\rho)}{M_{o,u}}\frac{m_{o,v}}{M_{o,v}v(P_\rho)},
\]
the Boundary Harnack Inequality \eqref{Eq:3:8bis} can be rewritten as
\[
\tilde\gm^{-1}\frac{u(P_\rho)}{v(P_\rho)}\le\frac{u(x,t)}{v(x,t)}\le\tilde\gm\frac{u(P_\rho)}{v(P_\rho)}
\]
where now $\tilde\gm$ depends not only on $N$, $p$, $M_{1,1}$, but also on $M_{o,u}/m_{o,u}$
and $M_{o,v}/m_{o,v}$.}
\end{remark}
%%%%%%%%%%%%%%%%%%%%%%%%%%%%%%%%%%%%%%%%%%%%%%%%
\begin{remark}
{\normalfont Note that {\bf (a)} implies that near a part of the lateral boundary, 
where a non-negative solution vanishes, the parabolic 
$p$--Laplace operator is uniformly elliptic. {Since we do not have an estimate
at the boundary of the type $|Du(x,t)|\ge\frac{u(x,t)}{d(x,\partial E)}$, {\bf (a)}
and {\bf (c)} hold only in a small neighbourhood of $S_T$, whose size depends on the solution,
as both $c_o$ and the oscillation of the gradient $Du$
depend on  the oscillation of $u$: this is precisely
the meaning of $\rho_o$. Moreover, as it will be clear from the proof, we give only a qualitative
dependence of $\be$ on the various quantities.}}
\end{remark}
{The proof relies on proper estimates from above and below, 
which were originally proved in \cite[\S~4]{DKV}
for solutions to the singular porous medium equations in $C^2$ domains
by building explicit barriers; they were
later extended in \cite{SV} to solutions to doubly
nonlinear singular equations in $C^{1,\al}$ domains. Unfortunately, it recently
turned out that there is a flaw in the argument, and, as pointed out in \cite{KMN}
in the context of $p$--laplacian type equations,
$C^{1,1}$ domains is the most general assumption one can have,
in order to build barriers; that this is the threshold below which one cannot go,
had already been shown in the elliptic context in \cite{Aik2}.}

We recast these estimates in the lemma below, 
in a form tailored to our purposes.
Indeed, the Hopf Principle and a weak version of the Boundary Harnack 
Inequality follow easily from these estimates.
Our improvement lies in the use of the Carleson estimates, that allow a more precise 
bound for $\frac{u(x,t)}{v(x,t)}$ in terms of $\frac{u\prho}{v\prho}$. The restriction 
to $\frac{2N}{N+1}<p<2 $ comes into play only in this last step.

Thus, let $\partial E$ be of class $C^{1,1}$ and $u$ be a non-negative, 
weak solution to \ref{Eq:1:1:o} in $E_T$, for
$1<p<2$. Assume that $u\leq M$ in $E_T$.
For $x\in \rn$, set $d(x)= \dist(x,\partial E)$, and for
$s>0$, let
$$
E^s=\{x\in E:\ \ \frac{s}{2}\leq d(x)\leq 2s\}.
$$

\begin{lemma}\label{Lm:3:1}
Let $\tau\in(0,T)$ and fix $x_o\in\partial E$. Assume that $u$ vanishes on $$\partial E\cap 
K_{2\rho}(x_o)\times(\tau,T).$$
For every $\nu>0$, there exist positive
constants $\gm_1,\gm_2$, and $0<\bar s<\frac12$, depending only on $N$, $p$, $\nu$, and
$M_{1,1}$, such that for all
$\tau+\nu M^{2-p}\rho^p<t<T$, and for all $x\in E\cap K_{2\bar s\rho}(x_o)$ with
$d(x)<\bar s\rho,$
\begin{equation*}%\label{Eq:4:19}
\gm_2\left(\frac{d(x)}\rho\right) \inf_{K_{2\rho}(x_o)\cap E^{\bar s\rho}\times(\tau,T)}u \le u(x,t)\le\gm_1 \left(\frac{d(x)}{\rho}\right)\sup_{E\cap K_{2\rho}(x_o)\times(\tau,T)} u.
\end{equation*}
\end{lemma}
Relying on the above lemma, the proof of Theorem~\ref{Thm:3:3} follows rather easily.
\smallskip

\subsubsection{\textit{Proof of Theorem~\ref{Thm:3:3}.}}
In order to simplify the proof, we assume $M_{1,1}\approx1$. Otherwise,
the anisotropies of the domains we have to take into account would make the
reading particularly burdensome.
\vskip.2truecm
\noi{\bf (a)} The Hopf Principle \eqref{Eq:3:8} is an easy consequence of
the linear growth estimate, the interior Harnack inequality, and the regularity 
up to the boundary of $Du$ (\cite[Chapters~IX and X]{dibe-sv}).
\vskip.2truecm
\noi{\bf (b)} To prove the Boundary Harnack Inequality, let $\pto\in S_T$,
$\rho\in(0,r_o)$, $x_\rho=(x_o',2M_{1,1}\rho)$, $P_\rho=(x_\rho,t_o)$.
With respect to Lemma~\ref{Lm:3:1}, let
\[
\tau=t_o-M^{2-p}\rho^p,\qquad T=t_o+M^{2-p}\rho^p,\qquad\nu=\frac12,
\]
and define
\begin{equation*}
V_{\frac12,\rho}(x_o,t_o)=\left\{x\in K_{\bar s\frac\rho4}(x_o)\cap E: d(x)<\bar s\frac\rho8\right\}
\times 
I(t_o,\rho,\frac12 M^{2-p}),
\end{equation*}
where $\bar s$ is the quantity claimed by Lemma~\ref{Lm:3:1} when $\nu=\frac12$.
Since both $u$ and $v$ satisfy the assumptions of Lemma~\ref{Lm:3:1}, 
for every $(x,t)\in V_{\frac12,\rho}(x_o,t_o)$, we can write
\begin{equation}\label{Eq:3:9}
\begin{aligned}
&\gm_2\,d(x)\,\mu_u(\bar s\rho)\le \rho u(x,t)\le\gm_1\,d(x)\,M_u(\rho),\\
&\gm_2\,d(x)\,\mu_v(\bar s\rho)\le \rho v(x,t)\le\gm_1\,d(x)\,M_v(\rho),
\end{aligned}
\end{equation}
where
\begin{align*}
&M_u(\rho)=\sup_{(E\cap K_{\frac\rho4}(x_o))\times I(t_o,\rho,M^{2-p})} u\\
&\mu_u(\bar s\rho)=\inf_{(E^{\bar s\frac\rho8}\cap K_{\frac\rho4}(x_o))\times I(t_o,\rho,M^{2-p})} u,
\end{align*}
and analogously for $v$; \eqref{Eq:3:9} yields
\begin{equation}  \label{Eq:3:10}
\frac{\gm_2}{\gm_1}\frac{\mu_u(\bar s\rho)}{M_v(\rho)}\le\frac{u(x,t)}{v(x,t)}
\le\frac{\gm_1}{\gm_2}\frac{M_u(\rho)}{\mu_v(\bar s\rho)},\qquad\forall (x,t)\in V_{\frac12,\rho}(x_o,t_o).
\end{equation}
Notice that \eqref{Eq:3:10} holds for every $1<p<2$. Restricting $p$ in the range $(\frac{2N}{N+1},2)$, allows us to apply Corollary~\ref{Cor:3:01} and Theorem~\ref{Thm:3:2}.

\noi Concerning $M_u(\rho)$, by Corollary~\ref{Cor:3:01} we have
\[
M_u(\rho)\le\gm\sup_{\tau\in I(t_o,\rho,2M^{2-p})}u(x_\rho,\tau).
\]
On the other hand, $\mu_u(\bar s\rho)$ is attained at some point 
$(x_*,t_*)\in K_{\frac\rho4}(x_o)\cap E^{\bar s\frac\rho8}\times I(t_o,\rho,M)$. By the elliptic Harnack inequality of 
Theorem~\ref{Thm:3:2}
\[
\gm_3 u(x_\rho,t_*)\le u(x_*,t_*)\le\gm_4 u(x_\rho,t_*),
\]
where $\gm_3$, $\gm_4$ depend on $N$, $p$, $M_{1,1}$ (as a matter of fact,
they depend on $\bar s$ too, 
but once $\nu$ is fixed, $\bar s$ depends only on these quantities), and therefore,
\[
\mu_u(\bar s\rho)\ge\gm_5\inf_{\tau\in I(t_o,\rho,2M^{2-p})}u(x_\rho,\tau),
\]
where, once more, $\gm_5$ depends only on $N$, $p$, $M_{1,1}$. 
Combining the previous estimates for $u$ and the analogous ones for $v$, yields
\begin{equation*}
\gm^{-1}\,\frac{\dsty\inf_{\tau\in I(t_o,\rho,2M^{2-p})}u(x_\rho,\tau)}{\dsty\sup_{\tau\in I(t_o,\rho,2M^{2-p})}v(x_\rho,\tau)}\le \frac{u(x,t)}{v(x,t)}\le\gm\,\frac {\dsty\sup_{\tau\in I(t_o,\rho,2M^{2-p})}u(x_\rho,\tau)}{\dsty\inf_{\tau\in I(t_o,\rho,2M^{2-p})}v(x_\rho,\tau)},
\end{equation*}
for all $(x,t)\in V_{\frac12,\rho}(x_o,t_o)$.
\vskip.2truecm
\noi{\bf (c)} We prove the H\"older continuity of
the quotient $u/v$ up to the boundary. First notice that, given $(x,
t)\in S_T\cap\Psi_{\frac{\rho_o}2,M}(x_o,t_o)$ and denoting by $z_x = x + r\nu_x$ a point along
the normal to $\partial E$ at $x$, we can write
\begin{equation*}
\frac{u(z_x, t)}{v(z_x, t)}=\frac{u(z_x,t) - u(x,t)}{v(z_x,t) - v(x,t)} =%
\frac{\frac{\partial u}{\partial \nu_x}(c_x,t)}{\frac{\partial v}{\partial \nu_x}(c_x,t)}
\end{equation*}
with a suitable $c_x$. Since $u$ vanishes on the $C^{1,1}$ boundary, the tangential
component of the gradient vanishes, and $\dsty\frac{\partial u}{\partial\nu_x}\equiv |Du|$.
Letting $r\to0$, by the Hopf principle and
the $C^{\al}$ regularity of $Du$ and $Dv$ (see \cite[Chapters~IX and X]{dibe-sv}), 
we infer that $u/v$ has a H\"older
continuous trace on $S_T\cap\overline{\Psi_{\frac{\rho_o}2,M}\pto}$.
Moreover, if $(y, s)\in S_T\cap\overline{\Psi_{\frac{\rho_o}2,M}\pto}$ and $(x,t)\in\Psi_{\frac{\rho_o}2,M}\pto$,
we have, once more by Hopf's principle and the $C^{\al}$ regularity of $Du$
and $Dv$,
\begin{equation}\label{Eq:3:10bis}
\Big|\frac{u(x,t)}{v(x,t)} - \frac{u(y,s)}{v(y,s)}\Big| \le Cd((x,
t),(y,s))^{\al},
\end{equation}
where $d$ denotes the parabolic distance. 
When both $(x,t)$ and $(y,s)$ belong to $	\Psi_{\frac{\rho_o}2,M}\pto$, strictly speaking we should
distinguish three cases:
\begin{description}
\item $d((x,t),(y,s))\simeq d((x,t),S_T)\simeq d((y,s),S_T)$; in such a case, the interior
H\"older continuity suffices to conclude;
\item $d((x,t),(y,s))\gg d((x,t),S_T)$, $d((x,t),(y,s))\gg d((y,s),S_T)$; in such a case, on the left-hand side of \eqref{Eq:3:10bis} we can add and subtract $u/v$ evaluated at the boundary;
\item $d((x,t),(y,s))\simeq d((x,t),S_T)$, $d((x,t),(y,s))\gg d((y,s),S_T)$; in such a case
the wanted result is a straightforward consequence of the triangle inequality.
\end{description}
In all three instances, the proof is quite standard and basically relies on \eqref{Eq:3:10bis}
and on interior H\"older
continuity of $u/v$. Notice that here we can use the
classical parabolic distance, as the $p$--Laplacian is now uniformly
elliptic, thanks to the Hopf Principle.\hfill\bbox

%%%%%%%%%%%%%%%%%%%%%%%%%%%%%%%%%%%%%%%%%%%%%%%%

\subsection{Proof of Theorem~\protect\ref{Thm:3:1}}\label{S:3:4}

Although the overall strategy of the proof of Theorem~\ref{Thm:3:1} 
is the same of the proof of
Theorem~\ref{Thm:2:1}, the singular case requires an adapted 
renormalization argument to control the vertical size of proper dyadic cylinders. 
We mainly concentrate on the differences between the two proofs.
The H\"older continuity up to the boundary employed in 
Corollary~\ref{Cor:2:1}  also holds in the super-critical singular case, and 
the change of variables introduced in \S~\ref{S:2:4:1} works for any $p>1$.
Therefore, proving Theorem~\ref{Thm:3:1} reduces to proving the following 
lemma.
\begin{lemma}\label{Lm:3:2}
Let $u$ be a weak solution 
to \eqref{Eq:1:1}--\eqref{Eq:1:2} such that $0<u\le M$ in $E_T$.
Assume that $\partial E$ is flat
with respect to $x_{\scriptscriptstyle N}$ and that 
\begin{equation}\label{Eq:3:11}
I(t_o,9\rho,M^{2-p})\subset(0,T].
\end{equation}
Suppose that $u$ vanishes on
$$\partial E\cap \{|x_i-x_{o,i}|<2\rho,|x_{\scriptscriptstyle N}|<8L\rho\}\times I(t_o,9\rho,M^{2-p}).$$
Then there exists a constant $\gm>0$, depending only on $\datap$ such that
\begin{equation*}%%\label{Eq:3:11:bis}
u(x,t)\le\gm \sup_{\tau\in I(t_o,\rho,2 \etarho^{2-p})}u(x_\rho,\tau)
\end{equation*}
$\dsty\forall(x,t)\in\left\{|x_i-x_{o,i}|<\frac\rho4,\ 0<x_{\scriptscriptstyle N}<2L\rho\right\}
\times I(t_o,\rho,\etarho^{2-p})$.
\end{lemma}
\vskip.2truecm
\noi{\it Proof} -  Since the previous flattening of the boundary does not 
affect the value of $M$, 
we still have %\begin{equation}\label{Eq:3:6}
$\dsty 0< u\le M$ %\end{equation}
in $E_T$ for some $M>0$. By \eqref{Eq:3:11}
we can define the set
\begin{equation*}
\Psi^*_{\rho}=K^*_{2\rho}(x_o)\times I(t_o,9\rho,\etarho^{2-p})\subset E_T,
\end{equation*}
where $K^*_{2\rho}(x_o)$ has been defined in \eqref{Eq:2:16},
and $\etarho$ is the first root of the equation
\[
\max_{\Psi^*_\rho}u=\etarho.
\]
We have already shown that at least one solution of such equation does exist.
Moreover, $\Psi^*_\rho\subset E_T$. 
The change of variable
\begin{equation*}
x\to\frac{x-x_o}{2L\rho},
\qquad t\to\frac1{\etarho^{2-p}}\frac{t-t_o}{\rho^p}
\end{equation*}
maps $K_{2\rho}^*(x_o)
\times I(t_o,9\rho,\etarho^{2-p})$ into $%
\tilde Q=\{|y_i|<\frac1L,\ 0<y_{\scriptscriptstyle N}<2\}\times(-{9^p},{9^p}]$, 
$x_\rho$ into $y_o=(0,\dots,0,1)$, $\tilde{K}_{\rho}(x_o)$ into $\tilde K_1=\{|y_i|<\frac1{2L}, 
|y_{\scriptscriptstyle N}|<1\}$, 
$K^*_{2\rho}(x_o)$ into $K^*_2(y_o)=\{|y_i|<\frac1L,\ 0<y_{\scriptscriptstyle N}<2\}$, and the
portion of the lateral boundary $S_T\cap \overline{\Psi^*_{\rho}}$ into 
$$\Xi=\left\{(y^{\prime},0):\ |y_i|<\frac1L\right\}\times\left(-{9^p},{9^p}\right].$$

\noi Denoting again by $(x,t)$ the transformed variables, and letting
$y_o=(0,\dots,0,1)$ in the remainder of the proof, the rescaled
function
\begin{equation*}
v(x,t)=\frac1{\etarho}\,u(2L\rho\, x+x_o,t_o+\etarho^{2-p}\rho^p t)
\end{equation*}
is a non-negative, weak solution to
\begin{equation*}
\partial_t v-\dvg\tilde{\bl{A}}(x,t,v,Dv)=0,
\end{equation*}
in $\tilde Q$, where it is easy to see that $\tilde{\bl{A}}$ satisfies
structure conditions analogous to \eqref{Eq:1:2}, and
\begin{equation*}
\forall\,(x,t)\in \tilde Q\qquad 0\le v\le1.
\end{equation*}
As in the proof of Theorem~\ref{Thm:2:1}, in order to simplify the notation and
without loss of generality, from here on we assume $L=1$.

Since the structure conditions have changed, we will denote
with $\eps_*$ and $\gm_*$ the corresponding constants of the Harnack
inequality, claimed by Theorem~\ref{Thm:3:2}. We will repeatedly apply %
\eqref{Eq:3:3}--\eqref{Eq:3:4}: due to all our assumptions, the only
condition we need to take into account each time, is that $\dsty %
K_{8R}(x_*)\subseteq E$, where $ x_*\in\{|x_i|<\frac12,\,0<x_{\scriptscriptstyle N}<1%
\}$, and $R$ depends on the context.

The following argument closely resembles the proof given in the elliptic context, the main 
difference being only the need to control the size of the time interval. 
On the other hand, we cannot simply repeat here the argument used in the degenerate case, as it heavily relies on the fact that $p>2$. 

Consider the set
%\begin{equation*}
$\dsty F_o=\{(x^{\prime},x_{\scriptscriptstyle N}):\ |x_i|<\frac12,\,x_{\scriptscriptstyle N}=1\}$,
%\end{equation*}
the points $P_{h^{\prime}}$, whose coordinates are given by
\begin{align*}
&x_i=\frac{h_i}8,\quad h_i=-3,-2,\dots,3,\quad i=1,2,\dots,N-1, \\
&x_{\scriptscriptstyle N}=1,
\end{align*}
and the $(N-1)$-dimensional cubes $\dsty K_{\frac14}(P_{h^{\prime}})\cap%
\{x_{\scriptscriptstyle N}=1\}$. Notice that

\begin{itemize}
\item the $(N-1)$-dimensional cubes give an equal-size decomposition of $F_o$;

\item due to their size and their distance from the boundary, $%
K_2(P_{h^{\prime}})\subseteq E$ and therefore we can apply \eqref{Eq:3:3}--\eqref{Eq:3:4}.
\end{itemize}
Consequently, for any
\begin{equation}  \label{Eq:3:13}
t_*\in[-\eps_* v^{2-p}(y_o,0),\eps_* v^{2-p}(y_o,0)],
\end{equation}
by \eqref{Eq:3:3} we have, 
\begin{equation*}
v(x,t_*)\le\gm_* v(P_{h^{\prime}},t_*)
\end{equation*}
for all $\dsty x\in K_{\frac14}(P_{h^{\prime}})\cap\{x_{\scriptscriptstyle N}=1\}$. On the other
hand, it is easy to see that at most
\begin{equation*}
v(P_{h^{\prime}},t_*)\le\gm_* v(y_o,t_*).
\end{equation*}
Therefore, for any $x\in F_o$,
\begin{equation*}
v(x,t_*)\le\gm_*^4 v(y_o,t_*).
\end{equation*}
Consider the slab
\begin{equation*}
S_o=\{(x^{\prime},x_{\scriptscriptstyle N}):\ |x_i|<\frac12,\,\frac78<x_{\scriptscriptstyle N}<1\}.
\end{equation*}
As we noticed above, for any $\bar x\in F_o$, $\dsty K_2(\bar x)\subseteq E$; 
therefore we can apply \eqref{Eq:3:3} at time level $t_*$, and conclude
that
\begin{equation*}
\forall\,\bar x\in F_o,\ \ \forall x\in K_{\frac14}(\bar x)\quad v(x,t_*)\le%
\gm_* v(\bar x,t_*).
\end{equation*}
Consequently, $\forall\,x\in S_o$
\begin{equation} \label{Eq:3:14}
v(x,t_*)\le\gm_*^5 v(y_o,t_*).
\end{equation}
Estimate \eqref{Eq:3:14} holds in particular for any $x\in F_1$, where
\begin{equation*}
F_1=\{(x^{\prime},x_{\scriptscriptstyle N}):\ |x_i|<\frac12,\, x_{\scriptscriptstyle N}=\frac78\}.
\end{equation*}
We can then iterate and conclude that
\begin{equation}  \label{Eq:3:15}
\begin{aligned} &\forall\,x\in S_k=\{(x',x_{\scriptscriptstyle N}):\
|x_i|<\frac12,\,(\frac78)^{k+1}<x_{\scriptscriptstyle N}<(\frac78)^k\}\\
&v(x,t_*)\le\gm_*^{k+5}v(y_o,t_*). \end{aligned}
\end{equation}
On the other hand, by \eqref{Eq:3:3}--\eqref{Eq:3:4}, for any $t_*$ as in 
\eqref{Eq:3:13},
\begin{equation}  \label{Eq:3:16}
v(y_o,t_*)\le\gm_* v(y_o,0).
\end{equation}
Combining \eqref{Eq:3:15}--\eqref{Eq:3:16} finally yields
\begin{equation*}
v(x,t)\le\gm_*^{k+6} v(y_o,0)
\end{equation*}
for all $\dsty(x,t)\in S_k\times[-\eps_* v^{2-p}(y_o,0),\eps_* v^{2-p}(y_o,0)%
]$.

\noi For any $\tau\in[-2,2]$ (and not just for $\tau=0$) we can repeat the
same argument and obtain
\begin{equation*}
v(x,t)\le\gm_*^{k+6}v(y_o,\tau)
\end{equation*}
for all $\dsty(x,t)\in S_k\times[\tau-\eps_* v^{2-p}(y_o,\tau),\tau+\eps_*
v^{2-p}(y_o,\tau)]$, for any $\tau\in[-2,2]$, provided that $v(y_o,\tau)>0$. %
\noi Therefore, setting
%\begin{equation*}
$\dsty {\cal M}_2=\sup_{\tau\in[-2,2]}v(y_o,\tau)$,
%\end{equation*}
yields
\begin{equation}  \label{Eq:3:17}
v(x,t)\le\gm_*^{k+6}{\cal M}_2
\end{equation}
for all $\dsty(x,t)\in S_k\times[-2,2]$. This plays the role of
\eqref{Eq:2:22bis} in the singular framework.
Now let
\begin{align*}
&Q_*=\left\{|x_i|<\frac14,\ 0<x_{\scriptscriptstyle N}<1,\ t\in[-1,1]\right\}, \\
&Q^*=\left\{|x_i|<\frac14,\ -1<x_{\scriptscriptstyle N}<1,\ t\in[-1,1]\right\}, \\
&\tilde Q_*=\left\{|x_i|<\frac12,\ 0<x_{\scriptscriptstyle N}<1,\ t\in[-2,2]\right\}, \\
&\tilde Q^*=\left\{|x_i|<\frac12,\ -1<x_{\scriptscriptstyle N}<1,\ t\in[-2,2]\right\}.
\end{align*}
Notice that we need to assume $|x_i|<\frac14$ in the definition of $Q_*$ and $Q^*$
for two closely connected reasons:
\begin{itemize}
\item the first point $P_1$ we are going to choose must lie in a properly small cylinder, 
so that the sequence $\{P_m\}$ is all contained in $\tilde Q_*$;
\item as in the degenerate case, we need to choose $k_o$ such 
that it depends only on the data.
\end{itemize}
We extend $u$ from $\tilde Q_*$ to $\tilde Q^*$
by odd reflection: by Lemma~\ref{Lm:2:4} $u$ is still a (signed) solution to %
\eqref{Eq:1:1}--\eqref{Eq:1:2}.

\noi Suppose there exists $\dsty P_1=(x_1,t_1)=(x_{1,1},x_{1,2},%
\dots,x_{\scriptscriptstyle 1,N},t_1) \in Q_*$Ä such that
\begin{equation}  \label{Eq:3:18}
u(P_1)\ge \gm_*^{k_o+6}{\cal M}_2;
\end{equation}
by \eqref{Eq:3:17}, it must be $0<x_{1,N}<(\frac78)^{k_o}$, $%
|x_i|<\frac14$, $t_1\in[-1,1]$. Consider the cylinder
\begin{equation*}
Q(P_1)=\left\{|x-x_1|<2(\frac78)^{k_o+1}\eps^{-l},\ \ t_1-(2(\frac78)^{k_o+1}%
\eps^{-l})^p<t\le t_1\right\},
\end{equation*}
where $\eps$ is the quantity claimed by Lemma~\ref{Lm:2:3}, and $l\in\nn$ is
to be fixed. Without loss of generality, we may assume $\eps<\frac78$. Now
let
\begin{equation*}
\om_o^{(1)}=\osc_{Q(P_1)}v.
\end{equation*}
We do not know the precise value of $\om_o^{(1)}$, but thanks to %
\eqref{Eq:3:18}, the construction of $u$ by odd reflection, and the
normalization of $v$, we surely have
\begin{equation}  \label{Eq:3:19}
2{\cal M}_2\gm_*^{k_o+6}\le\om_o^{(1)}\le 2.
\end{equation}
Provided $k_o$ is large enough, we have $Q(P_1)\subset \tilde Q^*$. Set $\sig%
_o^{(1)}=2(\frac78)^{k_o+1}\eps^{-l}$ and consider
\begin{equation*}
Q_o^{(1)}=K_{\sig_o^{(1)}}(x_1)\times(t_1-\theta_o^{(1)}(\sig%
_o^{(1)})^p,t_1), \qquad\text{where}\ \ \ \theta_o^{(1)}=\left(\frac{\om%
_o^{(1)}}A\right)^{2-p},
\end{equation*}
and $A$ is the quantity denoted by $c$ in Lemma~\ref{Lm:2:3}, which we assume
to be larger than $2$, without loss of generality. It is apparent that $%
Q_o^{(1)}\subset Q(P_1)\subset \tilde Q^*$. By Lemma~\ref{Lm:2:3}, we can
build a sequence
\begin{equation*}
\om_n^{(1)}=\dl\om_{n-1}^{(1)},\quad \theta_n^{(1)}=\left(\frac{\om_n^{(1)}}%
A\right)^{2-p},\quad
\sig_n^{(1)}=\eps\sig_{n-1}^{(1)},\quad Q_n^{(1)}=Q_{\sig%
_n^{(1)}}(\theta_n^{(1)}),
\end{equation*}
for all non-negative integers $n$. Such a sequence satisfies
\begin{equation*}
Q_{n+1}^{(1)}\subset Q_n^{(1)},\qquad \osc_{Q_n^{(1)}}v\le\om_n^{(1)}.
\end{equation*}
By iteration
\begin{equation*}
\osc_{Q_n^{(1)}}v\le\dl^n\om_o^{(1)}=\dl^n\osc_{Q(P_1)}v,
\quad\Rightarrow\quad\osc_{Q(P_1)}v\ge\frac1{\dl^n}\osc_{Q_n^{(1)}}v.
\end{equation*}
If we now choose $n=l$, and $l$ such that $\dl^{-l}>\gm_*^{5}$, by the
choice of $\sig_o^{(1)}$ we conclude that
\begin{equation*}
\om_o^{(1)}\ge2\gm_*^{k_o+11}{\cal M}_2
\end{equation*}
and this obviously improves the previous lower bound given by \eqref{Eq:3:19}.
As $v$ has been built by odd reflection, we conclude there must exist
$P(x_2,t_2)=(x_{2,1},x_{2,2},\dots,x_{2,N},t_2)\in Q(P_1)$ such that
\begin{equation*}
v(P_2)\ge \gm_*^{k_o+11}{\cal M}_2.
\end{equation*}
As before, by \eqref{Eq:3:17}, we have $0<x_{2,N}<(\frac78)^{k_o+10}$, $%
|x_i|<\frac14$, and also $t_1-(2(\frac78)^{k_o+1}\eps^{-l})^p<t_2<t_1$. 
Set
\begin{equation*}
Q(P_2)=\left\{|x-x_2|<2(\frac78)^{k_o+10}\eps^{-l},\ \
t_2-(2(\frac78)^{k_o+10}\eps^{-l})^p<t<t_2\right\}.
\end{equation*}
Once more, provided $k_o$ is large enough, we can assume that $Q(P_2)\subset
\tilde Q^*$. From now on, we proceed as in \S ~\ref{S:2:4:3}. By induction,
we obtain a sequence $\{P_m=(x_m,t_m)\}$, such that
\begin{equation}  \label{Eq:3:20}
v(P_m)\ge \gm_*^{k_o+6+5(m-1)}{\cal M}_2,
\end{equation}
and
\begin{equation*}
0<x_{m,N}<(\frac78)^{k_o+1+5(m-1)}.
\end{equation*}
Provided we choose $k_o$ large enough, the sequence $\{P_m\}$ is all
contained in the fixed cylinder $\tilde Q_*$: since $0<v<1$ and ${\cal M}_2\in(0,1]$
is a fixed quantity, \eqref{Eq:3:20} eventually leads to a contradiction.
Therefore, there exists $\tilde\gm$ that depends only on the data, such that
\begin{equation}  \label{Eq:3:21}
\forall(x,t)\in Q_*\qquad v(x,t)\le\tilde\gm {\cal M}_2
\end{equation}

\noi We now switch back to the original variables.
We conclude that there exists a constant $\tilde\gm$, depending only on the data, such that
\begin{equation*}
\forall(x,t)\in\left\{|x_i-x_{o,i}|<\frac\rho4,\,
0<x_{\scriptscriptstyle N}<2L\rho\right\}\times I(t_o,\rho,\etarho^{2-p})
\end{equation*}
\begin{equation}
u(x,t)\le \tilde\gm \sup_{\tau \in I(t_o,\rho,2\etarho^{2-p})} u(x_\rho,\tau).\tag*{$\blacksquare$}
\end{equation}
\subsubsection{Proof of Corollary~\ref{Cor:3:01}}
A close inspection of the previous proof shows that all the arguments continue to hold true, 
if we substitute the qualitative parameter $\etarho$ directly with $M$.\hfill\bbox
%%%%%%%%%%%%%%%%%%%%%%%%%%%%%%%%%%%%%%%%%%%%%%%
\subsection{Proof of Corollary~\ref{Cor:3:1}}
Taking into account the notation of \S~\ref{S:3:4}, it is enough to prove the following result.
\begin{lemma}
Let $u$ be a weak solution 
to \eqref{Eq:1:1}--\eqref{Eq:1:2} such that $0<u\le M$ in $E_T$.
Assume that $\partial E$ is flat
with respect to $x_{\scriptscriptstyle N}$ and that $\dsty I(t_o,9\rho,M^{2-p})\subset(0,T]$.
Suppose that $u$ vanishes on
$$\partial E\cap \{|x_i-x_{o,i}|<2\rho,|x_{\scriptscriptstyle N}|<8L\rho\}\times I(t_o,9\rho,M^{2-p}).$$
Then there exists a constant $\gm>0$, depending on $\datap$, and $\frac M{m_o}$, 
such that
\begin{equation*}
u(x,t)\le\gm\,u(P_{\rho_k})
\end{equation*}
$\dsty\forall(x,t)\in\Psi_{\rho_k,M}\pto$, for all $k=0,1,2,\dots$.
\end{lemma}
\noi{\it Proof} - 
We go back to the proof of Lemma~\ref{Lm:3:2},
at \eqref{Eq:3:21} with $\etarho$ substituted by $M$. 
Let $(x_o,t_o)\in S_T$ and assume $\dist(x,\partial E)=\rho$ 
(we are now taking $L=\frac12$, without loss of generality, in order to simplify the notation), 
and $I(t_o,9\rho,M^{2-p})\subset(0,T]$. We let
\begin{equation*}
\begin{aligned}
&I_k=(t_o,\sig_k,2M^{2-p}),\qquad M_k=\sup_{\tau\in I_k}u(x_{\rho_k},\tau),\\
&m_k=\inf_{\tau\in I_k}u(x_k,\tau),\qquad u_k=u(x_{\rho_k},t_o).
\end{aligned}
\end{equation*}
By the weak form of the Carleson estimate given in Corollary~\ref{Cor:3:01}, 
$\forall(x,t)\in\Psi_{\rho,M}\pto$
\begin{equation*}
u(x,t)\le\gm\sup_{\tau\in I_o}u(x_\rho,\tau),
\end{equation*}
which implies
\begin{equation*}
u(x,t)\le\gm\frac{M_o}{m_o}\,u\prho\le\gm\frac M{m_o}u\prho.
\end{equation*}
Analogously, working in a smaller box, $\forall(x,t)\in\Psi_{\rho_k,M}\pto$,
\begin{equation*}
u(x,t)\le\gm\sup_{\tau\in I_k}u(x_{\rho_k},\tau)\le\gm\frac{M_k}{m_k}u(P_{\rho_k}).
\end{equation*}
%Notice that here the height of $\Psi_{\rho_k}$ depends on $u\prho$, 
%and not on $u(P_{\rho_k})$.
%This can always be done without loss of generality, by possibly stretching $\Psi^*_\rho$,
%where $\eta_\rho$ is defined.
To prove the lemma, we will show by induction that 
\begin{equation} \label{Eq:3:23}
	\frac{M_k}{m_k} \leq \gamma^{2(N_0+1)}, \text{ and } N_k \leq N_{0}
\end{equation}
where $N_0 = 2\left(\frac M{m_o}\right)^{2-p}\frac{8^p}{\bar\eps}$, and 
$N_k = 2 \left(\frac M{\gamma^k m_k}\right)^{2-p}\frac{8^p}{\bar \epsilon} $ for any $k=1,2,\dots$.

Let us now consider $k = 1$.
In order to cover the segment $\{x_\rho\}\times[t_o,t_o+2M^{2-p}\rho^p]$ 
(and the same can be said
for $\{x_\rho\}\times[t_o-2M^{2-p}\rho^p,t_o]$), we need at most $N_o$ steps, 
where $N_o$ is given by
\begin{equation*}
t_o+\bar\eps N_o m_o^{2-p}\left(\frac\rho8\right)^p=t_o+2M^{2-p}\rho^p,
\end{equation*}
which yields
\begin{equation*}
N_o=2\left(\frac M{m_o}\right)^{2-p}\frac{8^p}{\bar\eps}.
\end{equation*}
Without loss of generality, we can assume that $N_o\in\nn$, possibly
by a slight modification of $\bar\eps$.

Now consider $x_{\rho_1}$: by the elliptic Harnack inequality, $\forall\,t\in I_1$
we have
\begin{equation*}
\begin{aligned}
&u(x_{\rho_1},t)\ge\gm^{-1}u(x_\rho,t)\quad\Rightarrow\quad m_1\ge\gm^{-1}m_o,\quad\Rightarrow\quad\frac1{\gm m_1}\le\frac1{m_o}\\
&u(P_{\rho_1})\ge\gm^{-1}u\prho,\quad\Rightarrow\quad u(P_{\rho_1})\ge\gm^{-1}m_o.
\end{aligned}
\end{equation*}
Taking once more into account the Harnack inequality, we have
\begin{equation*}
\begin{aligned}
&M_1\le\sup_{t\in I_o}u(x_{\rho_1},t)\le\gm^{N_o+1}u\prho\\
&m_1\ge\inf_{t\in I_o}u(x_{\rho_1},t)\ge\gm^{-(N_o+1)}u\prho,
\end{aligned}
\end{equation*}
which yields
\begin{equation*}
\frac{M_1}{m_1}\le\gm^{2(N_o+1)}.
\end{equation*}
Moreover, to cover $\dsty\{x_{\rho_1}\}\times\left[t_o,t_o+2M^{2-p}\sig_1^p\right]$,
we need to take at most $N_1$ steps, where $N_1$ is given by
\begin{equation*}
t_o+\bar\eps N_1 m_1^{2-p}\left(\frac7{8^2}\rho\right)^p=t_o+2M^{2-p}\left(\frac7{8\gm^{\frac{2-p}p}}\rho\right)^p,
\end{equation*}
which yields
\begin{equation*}
N_1=2\left(\frac M{\gm m_1}\right)^{2-p}\frac{8^p}{\bar\eps}\le2\left(\frac M{m_o}\right)^{2-p}\frac{8^p}{\bar\eps},
\end{equation*}
and therefore $N_1\le N_o$.

Let us now assume that $\frac{M_k}{m_k} \leq \gamma^{2(N_0+1)}$, and $N_k \leq N_0$ holds for some $k = 1,2,\ldots$.

\noi Now consider $x_{\rho_{k+1}}\in E$: by the elliptic Harnack inequality, $\forall\,t\in I_{k+1}$
we have
\begin{equation*}
\begin{aligned}
&u(x_{\rho_{k+1}},t)\ge\gm^{-1}u(x_{\rho_k},t)\quad\Rightarrow\quad m_{k+1}\ge\gm^{-1}m_k,\\
&u(P_{\rho_{k+1}})\ge\gm^{-1}u(P_{\rho_k}),\quad\Rightarrow\quad u(P_{\rho_{k+1}})\ge\gm^{-1}m_k.
\end{aligned}
\end{equation*}
Taking once more into account the Harnack inequality, we have
\begin{equation*}
\begin{aligned}
&M_{k+1}\le\sup_{t\in I_k}u(x_{\rho_{k+1}},t)\le\gm^{N_k+1}u(P_{\rho_k})\\
&m_{k+1}\ge\inf_{t\in I_k}u(x_{\rho_{k+1}},t)\ge\gm^{-(N_k+1)}u(P_{\rho_k}),
\end{aligned}
\end{equation*}
which yields
\begin{equation*}
\frac{M_{k+1}}{m_{k+1}}\le\gm^{2(N_k+1)}\le\gm^{2(N_o+1)},
\end{equation*}
and the ratio has not grown.

\noi Moreover, to cover $\dsty\{x_{\rho_{k+1}}\}\times\left[t_o,t_o+2{M}^{2-p}\sig_{k+1}^p\right]$,
we need to take at most $N_{k+1}$ steps, where $N_{k+1}$ is given by
\begin{equation*}
t_o+\bar\eps N_{k+1} m_{k+1}^{2-p}\left(\frac{\rho_k}{8}\rho\right)^p=t_o+2M^{2-p} (\sigma_k)^p,
\end{equation*}
which yields
\begin{equation*}
N_{k+1}=2\left(\frac M{\gm^{k+1} m_{k+1}}\right)^{2-p}\frac{8^p}{\bar\eps}\le N_k.
\end{equation*}

\noi By the induction principle we now have \eqref{Eq:3:23}.
Hence
$$\frac{M_k}{m_k}\le\gm^{2(N_o+1)}\qquad\forall k=1,2,\dots,$$
where $N_o$ depends only on $\datap$, and $\frac M{m_o}$.
\hfill\bbox
%%%%%%%%%%%%%%%%%%%%%%%%%%%%%%%%%%%%%%%%%%%%%%%
\section{The Porous Medium Equation}\label{S:4}
Consider quasi-linear, parabolic partial differential equations 
of the form
\begin{equation}\label{Eq:4:1}
u_t-\dvg\bl{A}(x,t,u, Du) = 0\quad
\text{ weakly in }\> E_T,
\end{equation}
where the function $\bl{A}:E_T\times\rr^{N+1}\to\rn$ 
is only assumed to be measurable
and subject to the structure conditions 
\begin{equation}\label{Eq:4:2}
\left\{
\begin{array}{l}
\bl{A}(x,t,u,Du)\cdot Du\ge mC_o |u|^{m-1}|Du|^2\\
|\bl{A}(x,t,u,Du)|\le mC_1|u|^{m-1}|Du|
\end{array}\right .\quad \text{ a.e. }\> 
(x,t)\in E_T,
\end{equation}
where $C_o$ and $C_1$ are given positive constants, and 
$m>0$. 
The prototype of such a class of parabolic 
equations is
\begin{equation}\label{Eq:4:1:o}
u_t-\Delta (|u|^{m-1} u) =0 
\quad\text{ weakly in }\> E_T.\tag*{(4.1)${}_o$}
\end{equation}
For simplicity, we limit ourselves to the definition of non-negative solutions: a non-negative function 
\begin{equation*}%\label{Eq:4:4}
\begin{array}{l}
u\in C\big([0,T];L^2(E)\big)\ \text{ with }\
u^{\frac{m+1}2}\in L^2\big(0,T; W^{1,2}(E)\big)\ \ \text{if}\ \ m>1,\\
 \\
u\in C\big([0,T];L^{2}(E)\big)\ \text{ with }\
u^m\in L^2\big(0,T; W^{1,2}(E)\big)\ \ \text{if}\ \ 0<m<1
\end{array}
\end{equation*}
is a weak sub(super)-solution to (\ref{Eq:4:1})--(\ref{Eq:4:2}) 
if for every sub-interval 
$[t_1,t_2]\subset (0,T]$
\begin{equation}\label{Eq:4:3}
\int_E u\vp dx\bigg|_{t_1}^{t_2}+\int_{t_1}^{t_2}\int_E
\big[-u\vp_t+\bl{A}(x,t,u,Du)\cdot D\vp\big]dxdt\le(\ge)0
\end{equation}
for all non-negative testing functions
\begin{equation*}
\vp\in W^{1,2}\big(0,T;L^2(E)\big)\cap 
L^2\big(0,T;W_o^{1,2}(E)\big).
\end{equation*}
This guarantees that all the integrals in (\ref{Eq:4:3}) 
are convergent. 

In \eqref{Eq:4:3}, the symbol $Du$ has to be understood in the sense of the following definition
\begin{equation*}
\begin{array}{l}
Du=\tfrac{2}{m+1}\mathbf 1_{\{ u>0\}} u^\frac{1-m}{2}Du^\frac{m+1}{2}\ \ \text{if}\ \ m>1,\\
 \\
Du=\tfrac{1}{m}\mathbf 1_{\{ u>0\}} u^{1-m}Du^{m}\ \ \text{if}\ \ 0<m<1.
\end{array}
\end{equation*}

The parameters $\datam$ are the data; the partial differential equation 
(\ref{Eq:4:1})--(\ref{Eq:4:2}) is degenerate when $m>1$, and singular when $0<m<1$,
since its modulus of ellipticity $|u|^{m-1}$ tends to $0$ or to $+\infty$, respectively, as $|u|\to0$. In the latter case, we further distinguish between singular super-critical range (when $\frac{(N-2)_+}N<m<1$), and singular critical and sub-critical range (when $0<m\le\frac{(N-2)_+}N$).

We are interested in solutions to \eqref{Eq:4:1}--\eqref{Eq:4:2} \emph{continuously}
vanishing on some distinguished part of the lateral boundary $S_T$ of
a space-time cylinder, and we aim at extending to such solutions 
the results proved in the previous sections for solutions to \eqref{Eq:1:1}--\eqref{Eq:1:2}. 
We will not give the full proofs, as most of the arguments
can be reproduced almost verbatim. We will limit ourselves to discuss the
changes, that need to be done.  

When dealing with this kind of problems, the cylinders to be considered are the ones 
already defined in \S~\ref{S:1}, provided the height scales as $\theta\rho^2$,
with $\theta=\left(\frac c{u\pto}\right)^{m-1}$, and $c$ a proper positive parameter.

Intrinsic Harnack inequalities, of the same kind as considered in 
Theorems~\ref{Thm:2:3} and \ref{Thm:3:2}, can be proved respectively for $m>1$, 
and $\frac{(N-2)_+}N<m<1$. For the exact statements we refer to
 \cite{DBGV-acta, DBGV-sing07,DBGV-mono}. For a general introduction to the porous medium equation, see \cite{DaKe,V}.
 %%%%%%%%%%%%%%%%%%%%%%%%%%%%%%%%%%%%%%%%%%%%%%%
\subsection{The Degenerate Case $m>1$}
For a Lipschitz cylinder 
$E_T$ with Lipschitz constant $L$, let $\pto$, $\rho$,  $P_\rho$ be as in \S~\ref{S:2:1},
assume $u\prho>0$, and set
\begin{equation*}
\begin{aligned}
&\Psi^-_\rho\pto\\
&=E_T\cap\{(x,t):\,|x_i-x_{o,i}|<\frac\rho4,\ |x_{\scriptscriptstyle N}|<2L\rho,
\ t\in(t_o-\frac{\al+\be}2\theta\rho^2,t_o-\be\theta\rho^2]\},
\end{aligned}
\end{equation*}
where $\theta=[\frac c{u\prho}]^{m-1}$, 
with $c$ given by the Harnack inequality for $m>1$ (see \cite{DBGV-mono},
Chapter~5), and $\al>\be$ are two positive parameters.
\begin{theorem}\label{Thm:4:1}
\emph{(Carleson Estimate, $m>1$)}
Let $u$ be a non-negative, weak solution to \eqref{Eq:4:1}--\eqref{Eq:4:2} in $E_T$.
Assume that
$$(t_o-\theta(4\rho)^2,t_o+\theta(4\rho)^2]\subset(0,T]$$ and that $u$ vanishes 
continuously on
$$\partial E\cap \{|x_i-x_{o,i}|<2\rho,\,|x_{\scriptscriptstyle N}|<8L\rho\}\times(t_o-\theta(4\rho)^2,t_o+\theta(4\rho)^2).$$
Then there exist two positive parameters $\al>\be$, and a constant $\gm>0$, 
depending only on
 $\datam$ and $L$, such that
\begin{equation}\label{Eq:4:5:o}
u(x,t)\le\gm\, u\prho\qquad\text{for every}\  \  (x,t)\in\Psi^-_{\rho}(x_{o},t_{o}).
\end{equation}
\end{theorem}
\vskip.2truecm
\noi{\it Proof} - As in the proof of Theorem~\ref{Thm:2:1}, the
flattening of the boundary leaves the structure of the equation unchanged:
relabeling the variables as before, we end up with
\begin{equation*}
\left\{
\begin{array}{l}
\bl{A}(x,t,u,Du)\cdot Du\ge mC^*_o |u|^{m-1}|Du|^2\\
|\bl{A}(x,t,u,Du)|\le mC^*_1|u|^{m-1}|Du|
\end{array}\right .\quad \text{ a.e. }\> 
(x,t)\in E_T
\end{equation*}
where {$C^*_o$ and $C^*_1$ are positive constants that depend only on $C_o$, $C_1$,
and the Lipschitz constant $L$}.
The rest of the proof proceeds as in Theorem~\ref{Thm:2:1}. \hfill\bbox
{\begin{remark}
{\normalfont As for the $p$--Laplacian, let us point out that also for the prototype 
equation \ref{Eq:4:1:o}, estimate \eqref{Eq:4:5:o}
could be extended from Lipschitz cylinders to a wider class of cylinders $E_T$, whose cross section $E$ is an {N.T.A. domain}.}
\end{remark}}
\vskip.2truecm 
\noi Weak solutions to \eqref{Eq:4:1} with zero
Dirichlet boundary conditions on a Lipschitz domain are H\"older continuous
up to the boundary (see, for example, \cite[Chapter~III, Appendix]{dibe-sv}). Combining this result with the previous Carleson estimate,
yields a quantitative estimate on the decay of $u$ at the boundary,
invariant by the intrinsic rescaling
\[
x=x_o+\rho y,\qquad t=t_o+\frac{\rho^2}{u\prho^{m-1}}\tau.
\]
\begin{corollary}\label{Cor:4:1}
Under the same assumptions of Theorem~\ref{Thm:4:1},
we have
\[
u(x,t)\le\gm\left(\frac{\dist(x,\partial E)}{\rho}\right)^{\mu}\,u\prho,
\]
for every $(x,t)\in\Psi^-_{\frac\rho2}\pto$, where $\gm>0$ and $\mu\in(0,1)$ depend only on 
$\datam$ and $L$.
\end{corollary}
\noi If we restrict our attention to solutions to the model equation \ref%
{Eq:4:1:o} and to $C^2$ cylinders, the result of Corollary~\ref{Cor:4:1} can
be strengthened. 
%%%%%%%%%%%%%%%%%%%%%%%%%%%%%%%%%%%%%%%%%%%%%%%%
\begin{theorem}%%\label{Thm:4:2}
\emph{(A Sharper Decay)} Let $E_T$ be a $C^2$ cylin\-der,
and $u$ a non-negative,
weak solution to \ref{Eq:4:1:o} in $E_T$.
Let the other assumptions of Theorem~\ref{Thm:4:1} hold.
Then there exist two positive parameters $\al>\be$, and a constant $\gm>0$, 
depending only on $m$, $N$, and
the $C^2$--constant $M_{2}$ of $E$,
such that for every 
$$(x,t)\in E_T\cap\left\{|x_i-x_{o,i}|<\frac{\rho}4,\ 
0<x_{\scriptscriptstyle N}< 2M_2\rho\right\}\times\left(t_o-\frac{\al+3\be}4\theta\rho^2,t_o-\be\theta\rho^2\right]$$
\begin{equation}\label{Eq:4:5}
0\le\,u(x,t)\le\gm\, \left(\frac{\dist(x,\partial E)}{\rho}\right)^{\frac1m}u\prho,
\end{equation}
\end{theorem}
\noi{\it Proof} - In this context the barrier is
\begin{equation*}
\Theta_k(x,t)=Cu\prho(1-\eta_k(x,t))^{\frac1m},
\end{equation*}
where 
\begin{equation*}
\eta_k(x,t)=\exp(-k(|x-y|-1))\exp(u\prho^{m-1}(t-t_1)),
\end{equation*}
and the constant $C$ is chosen so that $u\le\Theta_k$ on the parabolic boundary
of ${\cal N}_k$. The remainder of the proof is as in Theorem~\ref{Thm:2:2}; 
see also \cite[Theorem~4.1]{DKV}.
\hfill\bbox
\vskip.2truecm
Notice that, in general, the bound below by zero in \eqref{Eq:4:5} cannot be
improved. Indeed, when $m>1$, 
two explicit solutions to \ref{Eq:4:1:o} in the half-space 
$\{x_{\scriptscriptstyle N}\ge0\}$, that vanish at $x_{\scriptscriptstyle N}=0$, are
given by
\begin{equation*} 
%\begin{array}{l}
u_1(x,t)=x_{\scriptscriptstyle N}^{\frac1m},\qquad
u_2(x,t)=\left(\frac{m-1}{m+1}\right)^{\frac1{m-1}}(T-t)^{-\frac1{m-1}}x_{\scriptscriptstyle N}^{\frac2{m-1}}.
\end{equation*}
{As in \eqref{Eq:2:6}, we can even have solutions that vanish on a set
of positive measure. For $\gm\in(0,1)$, consider, for example,
\begin{equation*}
u(x,t)=\left[\frac{m-1}m\,\gm\,(t+1)\left(\gm+\frac{x_{\scriptscriptstyle N}-2}{t+1}\right)_+\right]^{\frac1{m-1}},
\end{equation*}
which solves \ref{Eq:4:1:o} in $\{x_{\scriptscriptstyle N}>0\}\times(0,\frac2\gm-1]$.}

Finally, all the remarks discussed in \S~\ref{SS:2:2} about the obstruction to the
definition of a useful Harnack chain, hold for the porous medium equation as well, 
without any significant change.
%%%%%%%%%%%%%%%%%%%%%%%%%%%%%%%%%%%%%%%%%%%%
\subsection{The Singular Super-critical Case $\frac{(N-2)_+}{N}<m<1$}

Let $E_T$, $u$, $\pto$, $\rho$, $P_\rho$ be as in Theorem~\ref{Thm:4:1}, and set
\begin{equation*}
\tilde{I}(t_o,\rho,h)=(t_o- h \rho^2,t_o+ h \rho^2).
\end{equation*}
Moreover, let $u$ be a weak solution to \eqref{Eq:4:1}--\eqref{Eq:4:2}
such that 
\begin{equation}\label{Eq:4:6}
0<u\le M\qquad\text{ in }\ \ E_T,
\end{equation} 
and assume that
\begin{equation}\label{Eq:4:7}
\tilde I(t_o,9\rho,M^{1-m})\subset(0,T].
\end{equation}
Then we define
\begin{equation*}
\begin{aligned}
&\widetilde\Psi_{\rho}\pto=E_T\cap\left\{(x,t):\ |x_i-x_{o,i}|<2\rho,\ |x_{\scriptscriptstyle N}|<4L\rho,
t\in \tilde I(t_o,9\rho,\etarho^{1-m})\right\},\\
&\bar\Psi_{\rho}\pto=E_T\cap\left\{(x,t):\ |x_i-x_{o,i}|<\frac\rho4,\ |x_{\scriptscriptstyle N}|<2L\rho,
t\in \tilde I(t_o,\rho,\etarho^{1-m})\right\},
\end{aligned}
\end{equation*}
where $\etarho$ is the first root of the equation
\begin{equation}  \label{Eq:4:8}
\max_{\widetilde\Psi_{\rho}\pto}u=\etarho.
\end{equation}
Working as in \S~\ref{S:3}, it is immediate to conclude that at least one root of \eqref{Eq:4:8}
does exist. Moreover, by \eqref{Eq:4:7} $\widetilde\Psi_{\rho}\pto\subset E_T$. 
Finally, for $k=0,1,2,\dots$ let
\begin{equation*}
\begin{aligned}
&\rho_k=\left(\frac78\right)^k\rho,\qquad\sig_k=\frac{\rho_k}{\gm^{k\frac{1-m}2}},\\ 
&x_{\rho_k}=(x_o',2L\rho_k),\qquad P_{\rho_k}=(x_o',2L\rho_k,t_o)\\
&\Psi_{\rho,M}\pto\\
&\quad=E_T\cap\left\{(x,t):\ |x_i-x_{o,i}|<\frac\rho4,\ |x_{\scriptscriptstyle N}|<2L\rho,
t\in \tilde I(t_o,\rho,M^{1-m})\right\},\\
&\Psi_{\rho_k,M}\pto\\
&\quad=E_T\cap\left\{(x,t):|x_i-x_{o,i}|<{\frac{\rho_k}4},\
|x_{\scriptscriptstyle N}|<2L\rho_k, \ t\in \tilde I(t_o,\sig_k,M^{1-m})\right\},\\
&m_o=\inf_{\tau\in \tilde I(t_o,\rho,2M^{1-m})}u(x_\rho,\tau),\quad M_o=
\sup_{\tau\in \tilde I(t_o,\rho,2M^{1-m})}u(x_\rho,\tau).
\end{aligned}
\end{equation*}
In the following statement, we give both the weak and the strong form
of the Carleson estimates for solutions to singular, super-critical porous
medium equation.
\begin{theorem}\label{Thm:4:3}
\emph{(Carleson Estimate, $\frac{(N-2)_+}{N}<m<1$)}. 
Let $u$ be a weak solution to \eqref{Eq:4:1}--\eqref{Eq:4:2},
that satisfies 
\eqref{Eq:4:6}. Assume that \eqref{Eq:4:7} holds true and $u$ vanishes 
continuously on
$$\partial E\cap \{|x_i-x_{o,i}|<2\rho,|x_{\scriptscriptstyle N}|<8L\rho\}\times 
\tilde{I}(t_o,9\rho,M^{1-m}).$$
Then there exist constants $\gm>0$ and $\al\in(0,1)$, depending only on
 $\datam$ and $L$, such that
 \begin{equation*}
u(x,t)\le\gm\left(\frac{\dist(x,\partial E)}{\rho}\right)^{\al}
\times\sup_{\tau\in I(t_o,\rho,2\etarho^{1-m})}u(x_\rho,\tau),
\end{equation*}
for every $(x,t)\in\bar\Psi_{\rho}\pto$.
Moreover, for every $(x,t)\in\Psi_{\rho,M}\pto$.
\begin{equation*}
u(x,t)\le\gm\left(\frac{\dist(x,\partial E)}{\rho}\right)^{\al}
\times\sup_{\tau\in I(t_o,\rho,2M^{1-m})}u(x_\rho,\tau),
\end{equation*}
Finally, there exists a constant $\hat \gm$, depending on $\datam$, $L$, and $\frac M{m_o}$,
such that
\begin{equation*}
u(x,t)\le \hat \gm\, u(P_{\rho_k}),
\end{equation*}
for every $(x,t)\in\Psi_{\rho_k}\pto$, for all $k=0,1,2,\dots$.
\end{theorem}
%%%%%%%%%%%%%%%%%%%%%%%%%%%
Another difference with respect to the degenerate case 
appears when we consider $C^{1,1}$ cylinders and (mainly for simplicity) 
the prototype equation \ref{Eq:4:1:o}. In this case, indeed, weak solutions
vanishing on the lateral part enjoy a proper power-like behavior at the boundary. 
{As we pointed out in \S~\ref{SS:strongCarleson},
the following Lemma was originally proved in \cite[\S~4]{DKV}
for solutions in $C^2$ domains; it can be extended to solutions in $C^{1,1}$ 
domains, working as done in \cite[Lemma~3.1]{KMN}
for the parabolic $p$--Laplacian.}
The notation is the same as for Lemma~\ref{Lm:3:1}.
\begin{lemma}\label{Lm:4:1}
Let $\tau\in(0,T)$ and fix $x_o\in\partial E$. Assume that $u$ vanishes on 
$\dsty\partial E\cap K_{2\rho}(x_o)\times(\tau,T)$.
For every $\nu>0$, there exist positive
constants $\gm_1,\gm_2$, and $0<\bar s<\frac12$, depending only on $N$, $m$, $\nu$, and
$M_{1,1}$, such that for all
$\tau+\nu M^{1-m}\rho^2<t<T$, and for all $x\in E\cap K_{2\bar s\rho}(x_o)$ with
$d(x)<\bar s\rho,$
\begin{equation*}%\label{Eq:4:19}
\gm_2\left(\frac{d(x)}\rho\right) \inf_{(E^{\bar s\rho}\cap K_{2\rho}(x_o)) \times(\tau,T)}u \le [u(x,t)]^m\le\gm_1 \left(\frac{d(x)}{\rho}\right)\sup_{(E\cap K_{2\rho}(x_o))\times(\tau,T)} u.
\end{equation*}
\end{lemma}
%%%%%%%% 
The implications of Lemma~\ref{Lm:4:1} are expressed in the following result.
\begin{theorem}%%\label{Thm:4:4}
Let $\frac{(N-2)_+}{N}<m<1$. Assume $E_T$ is a $C^{1,1}$ cylinder,
and  $\pto$, $\rho$, $P_\rho$ be as in Theorem~\ref{Thm:4:1}. 
Let $u,\,v$ be two weak solutions to \ref{Eq:4:1:o} in $E_T$, 
satisfying the hypotheses of Theorem~\ref{Thm:4:3}, $0<u,v\le M$ in $E_T$.
Then there exist positive  constants $\bar s$, $\gm$, $\be$, $0<\be\le1$,
depending only on $N$,
$m$, and $M_{1,1}$, and $\rho_o$, $c_o>0$, depending also on the oscillation of $u$,
such that the following properties hold.
\begin{description}
\item[(a)] \emph{Hopf Principle:}
\begin{equation*}
|Du^m|\ge c_o\qquad\text{in}\qquad  \Psi_{\rho_o,M}(x_o,t_o).
\end{equation*}
\item[(b)] \emph{Boundary Harnack Inequality:}
\begin{equation}\label{Eq:4:9}
\gm^{-1}\,\frac{\dsty\inf_{\tau\in \tilde I(t_o,\rho,2M^{1-m})}u(x_\rho,\tau)}{\dsty\sup_{\tau\in \tilde I(t_o,\rho,2M^{1-m})}v(x_\rho,\tau)}\le \frac{u(x,t)}{v(x,t)}\le\gm\,\frac {\dsty\sup_{\tau\in \tilde I(t_o,\rho,2M^{1-m})}u(x_\rho,\tau)}{\dsty\inf_{\tau\in \tilde I(t_o,\rho,2M^{1-m})}v(x_\rho,\tau)},
\end{equation}
for all $(x,t)\in\{x\in K_{\bar s\frac\rho4}\cap E:\,d(x)<\bar s\frac\rho8\}\times \tilde I(t_o,\rho,\frac12M^{1-m})$ with $\rho<\rho_o$.
\end{description}
\end{theorem}
\begin{remark}
{\normalfont Proceeding as in Remark~\ref{Rmk:3:4},
the Boundary Harnack Inequality \eqref{Eq:4:9} can be rewritten as
\[
\gm_*^{-1}\frac{u(P_\rho)}{v(P_\rho)}\le\frac{u(x,t)}{v(x,t)}\le\gm_*\frac{u(P_\rho)}{v(P_\rho)},
\]
where $\gm_*$ depends not only on $N$, $p$, $M_{1,1}$, but also on $M_{o,u}/m_{o,u}$
and $M_{o,v}/m_{o,v}$.}
\end{remark}

{\begin{remark}
{\normalfont As in Theorem~\ref{Thm:3:3}, the Hopf Principle holds true only in a small neighbourhood of the boundary.}
\end{remark}}

{\begin{remark}
{\normalfont Since we do not have H\"older regularity estimates for the gradient of $u$, we cannot proceed as in Theorem~\ref{Thm:3:3}, to prove the H\"older continuity of the ratio.}
\end{remark}}
%%%%%%%%%%%%%%%%%%%%%%%%%%%%%%%%%%%%%%%%%%%

%%%%%%%%%%%%%%%%%%%%%%%%%%%%%%%%%%%%%%%%%%%%%%%%
\bye
\end{document}